\documentclass[final]{siamart190516}

\usepackage{graphicx}
\usepackage{amsmath}
\usepackage{amsfonts}
\usepackage{amssymb}
\usepackage{url}
\usepackage{float}
\usepackage{tabu}
\usepackage{multirow}
\usepackage{rotating}
\usepackage{mathtools}
\usepackage{arydshln}
\usepackage{cleveref}

\usepackage{subfig}
\usepackage{adjustbox}

\usepackage{pifont}
\newcommand{\tick}{\ding{51}}
\newcommand{\cross}{\ding{55}}

\usepackage{color}
\newsiamremark{remark}{Remark}
\renewcommand{\vec}{\mathbf}
\definecolor{darkpastelgreen}{rgb}{0.01, 0.75, 0.24}
\definecolor{burgundy}{rgb}{0.5, 0.0, 0.13}

\newcommand{\VD}[1]{{\color{black} #1}}

\newcommand{\TheTitle}{A comparison of coarse spaces for Helmholtz problems in the high frequency regime}
\newcommand{\TheAuthors}{N. Bootland, V. Dolean, P. Jolivet, P.-H. Tournier}

\headers{A comparison of coarse spaces for Helmholtz problems}{\TheAuthors}

\title{{\TheTitle}\thanks{The first two authors gratefully acknowledge support from the EPSRC grant EP/S004017/1. The fourth author gratefully acknowledges support from the EPSRC grant EP/R009821/1.}}%

\author{N. Bootland%
	\thanks{Department of Mathematics and Statistics, University of Strathclyde, Glasgow, UK (\email{niall.bootland@strath.ac.uk}).}%
	\and
	V. Dolean%
	\thanks{Department of Mathematics and Statistics, University of Strathclyde, Glasgow, UK and Laboratoire J.A.~Dieudonn\'e, CNRS, University C\^ote d'Azur, Nice, France (\email{work@victoritadolean.com}).}%
	\and
	P. Jolivet%
	\thanks{IRIT, University of Toulouse, France (\email{pierre.jolivet@enseeiht.fr}).}%
	\and
	P.-H. Tournier%
	\thanks{LJLL, Sorbonne University, France (\email{tournier@ljll.math.upmc.fr}).}%
}

\Crefname{ALC@unique}{Line}{Lines}

\ifpdf
\hypersetup{
	pdftitle={\TheTitle},
	pdfauthor={\TheAuthors}
}
\fi

\begin{document}

\maketitle
\begin{abstract}
	Solving time-harmonic wave propagation problems in the frequency domain and within heterogeneous media brings many mathematical and computational challenges, especially in the high frequency regime. We will focus here on computational challenges and try to identify the best algorithm and numerical strategy for a few well-known benchmark cases arising in applications. The aim is to cover, through numerical experimentation and consideration of the best implementation strategies, the main two-level domain decomposition methods developed in recent years for the Helmholtz equation. The theory for these methods is either out of reach with standard mathematical tools or does not cover all cases of practical interest. More precisely, we will focus on the comparison of three coarse spaces that yield two-level methods: the grid coarse space, DtN coarse space, and GenEO coarse space. We will show that they display different pros and cons, and properties depending on the problem and particular numerical setting.
\end{abstract}






\section{Introduction}
\label{sec:Intro}

This work is motivated by the computational challenges that typically arise in frequency domain simulations of wave propagation and scattering problems in heterogeneous media. Such problems appear in a broad range of engineering applications, including acoustics, electromagnetics, and seismic inversion.

The discretisation of models describing frequency domain wave problems utilising finite element methodology almost always results in large, indefinite, and ill-conditioned linear systems. These linear systems are difficult to solve using standard methods, particularly for high frequencies and in the presence of complex heterogeneities. In order to maintain accuracy, the number of grid points must grow as a function of the frequency in such a way that, for high frequency problems, the size of the linear systems to be solved becomes prohibitive for direct methods. In such a regime, carefully designed iterative methods are required. Here, we consider a two-level domain decomposition approach for the robust parallel solution of the linear systems.

To model the wave problem, we utilise the Helmholtz equation on a domain $\Omega \subset \mathbb{R}^{d}$, $d = 2, 3$, for the field $u(\boldsymbol{x}) \colon \Omega \rightarrow \mathbb{C}$ given by
\begin{subequations}
\label{HelmholtzSystem}
\begin{align}
\label{HelmholtzEquation}
-\Delta u - k^{2} u & = f & & \text{in } \Omega,\\
\label{HelmholtzBCs}
\mathcal{C}(u) & = 0 & & \text{on } \partial\Omega,
\end{align}
\end{subequations}
where $\mathcal{C}$ incorporates some appropriate boundary conditions, $k(\boldsymbol{x}) > 0$ is the wave number, and $f(\boldsymbol{x})$ is a suitable forcing function. A key parameter is the wave number $k$, which relates the angular frequency $\omega$ and the wave speed $c$ as $k = \omega/c$. The wave speed $c(\boldsymbol{x})$ depends on the position $\boldsymbol{x}$ in the media for heterogeneous problems. Since $k$ is proportional to the frequency, the high frequency regime constitutes the case of large $k$ and presents particular challenges for designing effective solvers.

The difficulty in designing a good solver for the Helmholtz equation is shown very clearly in the review papers \cite{Ernst:2012:NAM, Gander:2018:SIREV} where one can see that there are no straightforward extensions to state-of-the-art methods for symmetric positive definite problems that tackle the indefinite or non-self adjoint problem well. Nonetheless, for large problem sizes---the case when one discretises the Helmholtz equation accurately for high wave numbers---domain decomposition methods are a natural choice \cite{Dolean:15:DDM}. However, despite recent efforts and in view of the latest results obtained both at the theoretical \cite{Graham:2017:RRD,Graham:2018:DDI,Gong:2020:DDP} or numerical level \cite{Dolean:2020:IFD,Dolean:2020:LSF,8818089}, there is no established method outperforming all others in the case of the Helmholtz problem.

Domain decomposition methods are well suited to solve large systems of equations arising from discretisation of PDEs and are among the best-known strategies for many types of problem. However, classical versions fail to be effective and may diverge for wave propagation problems. Two key constituent parts require a more careful treatment: the transmission conditions used to transfer information between adjoining subdomains and the coarse space that allows for capturing of global behaviour and passing information between subdomains globally. In this work we consider overlapping Schwarz methods.

The use of different transmission conditions at the interfaces (artificial boundaries arising from the decomposition into subdomains) has been extensively studied over the past two decades and various works \cite{Chevalier:1998:SMO,Collino:1997:NIM,Gander:2001:OSH} show that these conditions can improve the convergence of Schwarz methods and preconditioners. However, good transmission conditions are not sufficient to ensure a robust behaviour with respect to heterogeneities in the problem to solve or when the number of subdomains increases. To tackle these difficulties, we need coarse information that is cheap to compute and immediately available to all subdomains.

The focus in this work is on appropriate coarse spaces. A coarse space is typically required to provide scalability with respect to the number of subdomains used. More recently, however, coarse spaces have been designed to provide robustness to model parameters, especially for large contrasts in complex heterogeneous problems. For example, the GenEO (Generalised Eigenproblems in the Overlap) coarse space has been successfully employed for the robust solution of highly heterogeneous elliptic problems \cite{Spillane:2014:ARC,Haferssas:2017:ADS}. For the Helmholtz equation, finding a suitable coarse space is not an easy task and, being an indefinite problem, choosing a larger coarse space need not improve performance \cite{Fish:2000:GBT}. In designing coarse spaces for Helmholtz problems, we might also wish to reduce the dependence of the domain decomposition method on the wave number $k$. A natural idea to capture global behaviour is to use plane waves as a basis for the coarse space but it is not clear that this is suitable for heterogeneous media. Plane waves were first used within the multigrid approach \cite{Brandt:1997:WRM} before later being used to build coarse spaces for domain decomposition methods. We can cite the example of FETI(-DP)-H methods \cite{Farhat:2000:ATL,Farhat:2005:FDP}, for instance, but they have also been used in other domain decomposition methods \cite{Kimn:2007:ROB}. Nonetheless, plane waves have mainly been employed for homogeneous problems and do not have a straightforward extension to the heterogeneous case.

\VD{Even if coarse space information needs to be global and available to all domains, coarse spaces can be built locally and} based on local functions. Spectral coarse spaces use basis vectors deriving from the solution of local eigenvalue problems associated with appropriate operators. Within the context of the Helmholtz equation, this is exemplified in the DtN coarse space \cite{Conen:2014:ACS}. Here, eigenproblems are formulated on subdomain interfaces based on a Dirichlet-to-Neumann (DtN) map, extending an approach for elliptic problems \cite{Nataf:2011:ACS,Dolean:2012:ATL}. In this work we consider two spectral coarse spaces, the DtN approach and a GenEO-type approach suited to the Helmholtz problem. We will also consider a grid coarse space approach which utilises the addition of absorption in the problem.

Our consideration of coarse spaces for Helmholtz problems in the high frequency regime provide the following main contributions of the paper:
\begin{itemize}
	\item We bring together and outline recent work on developing coarse spaces that can be used to enhance domain decomposition methods for the Helmholtz problem in heterogeneous media. The approaches considered are then implemented in a common software, namely FreeFEM.
	\item We discuss implementation details and practical aspects of the methods as well as contrasting the benefits and drawbacks.
	\item We provide extensive numerical results on several well-known benchmark problems in 2D and 3D and compare the different approaches within a variety of settings.
	\item Based on the results of our numerical tests and our understanding of the implementation aspects involved, we provide an outlook on scenarios where certain methods may be more, or less, favourable.
\end{itemize}

The outline for the remainder of this work is as follows. In \Cref{sec:HelmholtzProblem} we detail the boundary value problem considered, i.e., the heterogeneous Helmholtz problem and its discretisation by finite elements. In \Cref{sec:DomainDecomposition} we introduce the basic principles of domain decomposition methods and present two versions of the one-level method, namely RAS and ORAS. A second level, or coarse space, is usually added by deflation. The three different coarse space strategies, namely the grid coarse space, DtN coarse space, and GenEO coarse space, are introduced in \Cref{sec:CoarseSpaces}. Parallel implementation details of these methods are given in \Cref{sec:ImplementationDetails} and an extensive numerical study is provided in \Cref{sec:ComparativeNumericalStudies}. Conclusions are then given in \Cref{sec:Conclusions}.

\section{The heterogeneous Helmholtz problem}
\label{sec:HelmholtzProblem}

Our model problem consists of solving the interior Helmholtz equation \eqref{HelmholtzSystem}. To be concrete, we let $\Omega$ is a bounded polygonal domain and consider specific boundary conditions on $\Gamma = \partial\Omega$. We suppose $\Gamma$ is partitioned into a disjoint union $\Gamma = \Gamma_{D} \cup \Gamma_{N} \cup \Gamma_{R}$ where Dirichlet conditions are imposed on $\Gamma_{D}$, Neumann conditions on $\Gamma_{N}$ and a Robin condition on $\Gamma_{R}$. Namely, we wish to solve
\begin{subequations}
\label{HelmholtzSpecificSystem}
\begin{align}
\label{HelmholtzSpecificEquation}
-\Delta u - k^{2} u & = f & & \text{in } \Omega,\\
\label{HelmholtzSpecificDirichletBC}
u & = u_{\Gamma_{D}} & & \text{on } \Gamma_{D},\\
\label{HelmholtzSpecificNeumannBC}
\frac{\partial u}{\partial \boldsymbol{n}} & = 0 & & \text{on } \Gamma_{N},\\
\label{HelmholtzSpecificRobinBC}
\frac{\partial u}{\partial \boldsymbol{n}} + i k u & = 0 & & \text{on } \Gamma_{R},
\end{align}
\end{subequations}
where $u_{\Gamma_{D}}$ is known. The Robin condition in \eqref{HelmholtzSpecificRobinBC} is a standard first order approximation to the far field Sommerfeld radiation condition and, in essence, enables appropriate wave behaviour to be described in a bounded domain, allowing for incoming or outgoing waves along $\Omega_{R}$. We do not require that a problem instance includes all types of boundaries but note that if $\Gamma_{R} = \emptyset$ then the problem will be ill-posed for certain choices of $k$. Furthermore, when $\Gamma_{R} \neq \emptyset$ the resulting linear systems, while being complex symmetric, are not Hermitian and this will be important in our choice of iterative method. However, classical iterative methods on their own are not enough to be able to solve Helmholtz problems effectively \cite{Ernst:2012:NAM}. This is further amplified when applied to highly heterogeneous problems.

The heterogeneity in our model is present in the wave number $k(\boldsymbol{x}) > 0$, being given by ratio of the angular frequency $\omega$ and the wave speed $c(\boldsymbol{x})$ as $k = \omega/c$. We allow $k$ to have jumps across different media and otherwise vary within the domain $\Omega$ such that $k \in L^{\infty}(\Omega)$.

To discretise \eqref{HelmholtzSpecificSystem}, we use the finite element method. In order to prescribe the weak formulation, we let $V = \left\lbrace u \in H^{1}(\Omega) \colon u = u_{\Gamma_{D}}\text{ on } \Gamma_{D}\right\rbrace$ and, in a similar fashion, $V_{0} = \left\lbrace u \in H^{1}(\Omega) \colon u = 0\text{ on } \Gamma_{D}\right\rbrace$. The weak form of the problem is then to find $u \in V$ such that
\begin{align}
\label{WeakForm}
a(u,v) & = F(v) & & \forall \ v \in V_{0},
\end{align}
where
\begin{subequations}
\begin{align}
\label{WeakFormBilinearPart}
a(u,v) & = \int_{\Omega} \left( \nabla u \cdot \nabla \bar{v} - k^2 u \bar{v}\right) \, \mathrm{d}\boldsymbol{x} + \int_{\Gamma_{R}} i k u \bar{v} \, \mathrm{d}s,
\end{align}
and
\begin{align}
\label{WeakFormLinearPart}
F(v) & = \int_{\Omega} f \bar{v} \, \mathrm{d}\boldsymbol{x},
\end{align}
\end{subequations}
are the bilinear and linear parts, respectively. To discretise, we consider piecewise polynomial finite element approximation on a simplicial mesh $\mathcal{T}^{h}$ of $\Omega$ which has a characteristic element diameter $h$. Denoting the associated trial space $V^{h} \subset V$ and test space $V_{0}^{h} \subset V_{0}$, the discrete problem is to find $u_{h} \in V^{h}$ such that
\begin{align}
\label{DiscreteWeakForm}
a(u_{h},v_{h}) & = F(v_{h}) & & \forall \ v_{h} \in V_{0}^{h}.
\end{align}
Let $\left\lbrace\phi_{j}\right\rbrace_{j=1}^{n}$ be the nodal basis for $V_{0}$ and $\left\lbrace\phi_{j}\right\rbrace_{j=n+1}^{n+d}$ be the nodal basis for the Dirichlet boundary $\Gamma_{D}$, for which $\mathcal{T}^{h}$ is assumed to conform. Then we can rewrite \eqref{DiscreteWeakForm} as a (complex) linear system
\begin{align}
\label{LinearSystem}
A\mathbf{u} = \mathbf{f},
\end{align}
where the coefficient matrix $A \in \mathbb{C}^{n \times n}$ and right-hand side vector $\mathbf{f} \in \mathbb{C}^{n}$ are given by $A_{i,j} = a(\phi_{j},\phi_{i})$ and $\mathrm{f}_{i} = F(\phi_{i}) - \sum_{l=n+1}^{n+d} a(\phi_{l},\phi_{i}) \bar{\mathrm{u}}_{l-n}$ respectively, for \mbox{$i,j \in 1, 2, \ldots, n$}. Here $\bar{\mathrm{u}}_{j}$, for $j = 1, 2, \ldots, d$, are the known Dirichlet values along $\Gamma_{D}$ corresponding to $u_{\Gamma_{D}}$. We then seek the solution $\mathbf{u} \in \mathbb{C}^{n}$ of the (in general) complex symmetric indefinite system \eqref{LinearSystem} to give
\begin{align}
\label{DiscreteSolution}
u_{h}(\boldsymbol{x}) = \sum_{j=1}^{n} \mathrm{u}_{j} \phi_{j}(\boldsymbol{x}) + \sum_{l=n+1}^{n+d} \bar{\mathrm{u}}_{l-n} \phi_{l}(\boldsymbol{x}).
\end{align}

The wave nature of solutions to the Helmholtz equation requires a sufficiently fine mesh in order to obtain a good approximation to the true solution and this should be kept in mind when considering the choice of discretisation of the problem. In terms of increasing the wave number $k$, if one is to maintain the same level of accuracy of discrete solutions then the number of grid points must increase faster than $k$ increases, due to the pollution effect \cite{Babuska:1997:IPE}. This growth depends on the discretisation chosen. For instance, in the case of using a piecewise linear (P1) finite element approximation on simplicial elements of diameter $h$, then $k^{3} h^{2}$ must be bounded, requiring $h$ to shrink as $\mathcal{O}(k^{-3/2})$. For piecewise quadratic (P2) finite elements on simplicial meshes the criteria relaxes to require that $h$ decreases as $\mathcal{O}(k^{-5/4})$. For higher order finite elements, the requirement becomes less stringent but the interpolation properties when using such approximation spaces ultimately begin to degrade.

Due to these restrictions and the desire for faster simulation times, it is common that practitioners simply consider a fixed number of points per wavelength instead, resulting in $h$ decreasing as $\mathcal{O}(k^{-1})$. Given a fixed number, $n_\text{ppwl}$, of points per wavelength we ensure that $n_\text{ppwl} h \approx \lambda$, where the wavelength is given by $\lambda = 2 \pi k^{-1}$. A prevalent engineering practice is to use 10 points per wavelength and in some large real-world problems of interest, such as in imaging science, it may be adequate or necessary to insist on less resolution. In light of this, we consider both 5 and 10 points per wavelength scenarios and make use of standard P2 finite element approximation throughout this work.

\section{Domain decomposition}
\label{sec:DomainDecomposition}

We now give details of the overlapping domain decomposition approach that we will utilise. This will be applied as a preconditioner rather than a stand-alone iterative method. Our approach is based on a two-level version of the optimised restricted additive Schwarz (ORAS) method. To provide a domain decomposition, we first partition $\Omega$ into non-overlapping subdomains $\left\lbrace\Omega'_{s}\right\rbrace_{s=1}^{N}$ which are resolved by the mesh $\mathcal{T}^{h}$. A layer of adjoining mesh elements is then added to provide overlapping subdomains $\left\lbrace\Omega_{s}\right\rbrace_{s=1}^{N}$ through the extension
\begin{align}
\label{OverlappingSubdomains}
\Omega_{s} = \mathrm{Int}\left(\bigcup_{\mathrm{supp}(\phi_{j})\cap\Omega'_{s}\neq\emptyset} \mathrm{supp}(\phi_{j})\right),
\end{align}
where $\mathrm{Int}(\cdot)$ denotes the interior of a domain and $\mathrm{supp}(\cdot)$ the support of a function. Note that more than one layer of elements can be added in a recursive manner if subdomains with larger overlap are required.

Now that we have a domain decomposition, we can define the restriction to a given subdomain $\Omega_{s}$ as an operator from $V^{h}$ into $V^{h}(\Omega_{s}) = \left\lbrace v\vert_{\Omega_{s}} \colon v \in V^{h} \right\rbrace$, namely $\mathcal{R}_{s} \colon V^{h} \rightarrow V^{h}(\Omega_{s})$ where $\mathcal{R}_{s} v = v\vert_{\Omega_{s}}$. Let $R_{s} \in \mathbb{R}^{n_{s} \times n}$ be the matrix form of $\mathcal{R}_{s}$ where $n_{s}$ is the number of degrees of freedom in $\Omega_{s}$. Since our subdomains overlap, we also make use of a partition of unity having matrix form $D_{s} \in \mathbb{R}^{n_{s} \times n_{s}}$ which is diagonal and satisfies $\sum_{s=1}^{N} R_{s}^{T} D_{s} R_{s} = I$; this removes ``double counting'' in the additive Schwarz method. Note that $R_{s}^{T}$ acts as an extension by zero outside of $\Omega_{s}$.

We can now define the restricted additive Schwarz (RAS) preconditioner
\begin{align}
\label{RAS}
M_{\text{RAS}}^{-1} = \sum_{s=1}^{N} R_{s}^{T} D_{s} A_{s}^{-1} R_{s},
\end{align}
where $A_{s} = R_{s}AR_{s}^{T}$ is the local Dirichlet matrix on $\Omega_{s}$, that is, Dirichlet conditions are implicitly assumed on $\partial\Omega_{s} \setminus \partial\Omega$. This corresponds to using Dirichlet transmission conditions in the additive Schwarz method. Note that each contribution from the sum in \eqref{RAS} can be computed locally in parallel.

The RAS approach, using the classical choice of Dirichlet transmission conditions, applied as a stand-alone stationary iterative method need not converge since frequencies in the error smaller than the wave number $k$ are not diminished \cite{Dolean:15:DDM}. When used as a preconditioner, the iterative solver used will typically suffer from slow convergence and may stagnate. Further, the local Dirichlet problems involving $A_{s}$ are not necessarily well-posed as $k^2$ may be an eigenvalue of the corresponding Laplace problem, in which case $A_{s}$ is singular. While methods to handle such singular systems can be applied, a different approach, which provides a convergent stand-alone method, is to change the Dirichlet transmission conditions to Robin conditions. This results in the so-called optimised restricted additive Schwarz (ORAS) method given by
\begin{align}
\label{ORAS}
M_{\text{ORAS}}^{-1} = \sum_{s=1}^{N} R_{s}^{T} D_{s} \widehat{A}_{s}^{-1} R_{s},
\end{align}
where now $\widehat{A}_{s}$ is the discretisation of the local Robin problem
\begin{subequations}
\label{ORASLocalSystem}
\begin{align}
\label{ORASLocalHelmholtzEquation}
-\Delta w_{s} - k^{2} w_{s} & = f & & \text{in } \Omega_{s},\\
\label{ORASLocalProblemBCs}
\mathcal{C}(w_{s}) & = 0 & & \text{on } \partial\Omega_{s}\cap\partial\Omega,\\
\label{ORASLocalRobinBC}
\frac{\partial w_{s}}{\partial \boldsymbol{n}_{s}} + i k w_{s} & = 0 & & \text{on } \partial\Omega_{s}\setminus\partial\Omega,
\end{align}
\end{subequations}
with $\mathcal{C}$ representing the underlying problem boundary conditions on $\partial\Omega$. The local problem \eqref{ORASLocalSystem} will always have a unique solution $w_{s}$. We note that the use of Robin conditions is not the only remedial choice. Optimal transmission conditions can also be studied and are given using a Dirichlet-to-Neumann (DtN) map \cite{Dolean:15:DDM}. However, this results in requiring pseudodifferential operators and is somewhat less practical without further approximation. We do not follow the pursuit of more advanced transmission conditions in this work.

The above approaches, RAS and ORAS, are one-level methods: they rely only on local subdomain solves and the local transmission of data. Such domain decomposition methods do not scale as we increase the number of subdomains used; that is, their convergence behaviour depends on $N$. To achieve robustness with respect to $N$, a coarse space is typically included, giving a two-level method. The coarse space is represented by a collection of column vectors $Z$, having full column rank. A coarse space operator $E = Z^{\dagger} A Z$ is constructed as well as the coarse correction operator $Q = Z E^{-1} Z^{\dagger}$; note the similarity to terms in \eqref{RAS}. The inclusion of the coarse space can be done in a number of ways, the simplest being additively as
\begin{align}
\label{2LevelAdditivePreconditioner}
M_{\text{2-level}}^{-1} = M_{\text{1-level}}^{-1} + Q,
\end{align}
where $M_{\text{1-level}}^{-1}$ is the underlying one-level preconditioner used, such as \eqref{RAS} or \eqref{ORAS}. Hybrid approaches are often more effective and we shall consider the adapted deflation technique
\begin{align}
\label{2LevelAdaptiveDeflationPreconditioner}
M_{\text{2-level}}^{-1} = M_{\text{1-level}}^{-1}P + Q = M_{\text{1-level}}^{-1}(I-AQ) + Q,
\end{align}
where $P = I-AQ$ is a projection. The most crucial choice is that of $Z$, which provides the coarse space. This is what we shall now consider.

\section{Coarse spaces}
\label{sec:CoarseSpaces}

The construction of a suitable coarse space can be achieved in different ways. One natural approach is to utilise a coarse grid in order to approximate the global behaviour. Often, the slow convergence in the one-level method can be characterised by ``slow modes'' which should be incorporated into the coarse space. For instance, in the Poisson problem slow modes correspond to constant functions in the kernel of the Laplace operator in each subdomain and these are used to give the Nicolaides coarse space. On the other hand, for homogeneous elasticity problems the slow modes are rigid body motions.

The notion that certain modes are responsible for slow convergence and must be incorporated into the coarse space can be made more general, in particular in the framework of spectral coarse spaces. Such coarse spaces use spectral information from an appropriate eigenproblem to identify relevant modes that should feature in the coarse space. Before considering spectral coarse spaces we first outline the grid coarse space.

\begin{remark}[Notation]
Following on from \Cref{sec:DomainDecomposition}, for a variational problem giving rise to a system matrix $B$ we denote by $B_{s}$ the corresponding local Dirichlet matrix on $\Omega_{s}$. Where Robin conditions are used on internal subdomain interfaces (the artificial boundaries) the local problem matrix is denoted by $\widehat{B}_{s}$. Meanwhile, when Neumann conditions are used on such interfaces, the local matrix is denoted by $\widetilde{B}_{s}$.
\end{remark}

\subsection{The grid coarse space method}
\label{subsec:CoarseGrid}
The grid coarse space was first introduced in \cite{Graham:2017:ADS} for the absorptive Helmholtz problem and extended to incorporate impedance (or Robin) conditions in \cite{Graham:2018:DDI}. In this case the one-level method is based on the following formula
\begin{align}
\label{eq:ORASeps}
M^{-1}_{1,\varepsilon} = \sum_{s=1}^{N}  R_s^T D_s \widehat{A}_{s,\varepsilon}^{-1} R_s.
\end{align}
where matrices $\widehat{A}_{s,\varepsilon}$ stem from the discretisation of the following local Robin problems with absorption (given by the parameter $\varepsilon \neq 0$)
\begin{align*}
-\Delta u_s - (k^2+i\varepsilon) u_s &= f  & & \text{in } \Omega_s,\\
\mathcal{C}(u_{s}) &= 0 & & \text{on } \partial\Omega_{s}\cap\partial\Omega,\\
\frac{\partial u_s}{\partial \boldsymbol{n}_s} + i k u_s &= 0  & & \text{on } \partial\Omega_{s}\setminus\partial\Omega.
\end{align*}
In order to achieve weak dependence on the wave number $k$ and number of subdomains $N$, the two-level preconditioner can be written in a generic way as follows
\begin{equation}
\label{eq:BNN}
M^{-1}_{2,\varepsilon} =  M^{-1}_{1,\varepsilon} P + Z E^{-1} Z^{\dagger},
\end{equation}
where $M^{-1}_{1,\varepsilon}$ is the one-level preconditioner \eqref{eq:ORASeps}, $Z$ is a rectangular matrix with full column rank, $E = Z^{\dagger}\widehat{A}_\varepsilon Z$ is the so-called coarse grid matrix, $Q = Z E^{-1} Z^{\dagger}$ is the so-called coarse grid correction matrix, and $P = I - A_\varepsilon Q$.

Perhaps the most natural coarse space is the one based on a coarser mesh, which we call the ``grid coarse space''. Let us consider $\mathcal{T}^{H_{\text{coarse}}}$, a simplicial mesh of $\Omega$ with mesh diameter $H_{\text{coarse}}$, and $V^{H_{\text{coarse}}} \subset V$, the corresponding finite element space. Let $\mathcal I_0 \colon V^{H_{\text{coarse}}} \to V^h$ be the nodal interpolation operator and define $Z$ as the corresponding matrix. Then, in this case, $E = Z^{\dagger}\widehat{A}_\varepsilon Z$ is really the stiffness matrix of the problem (with absorption) discretised on the coarse mesh.
Related preconditioners without absorption are used in \cite{Kimn:2007:ROB}.

\subsection{The DtN coarse space}
\label{subsec:DtN}

The Dirichlet-to-Neumann (DtN) coarse space \cite{Nataf:2011:ACS,Conen:2014:ACS} is based on solving local eigenvalue problems on subdomain boundaries related to the DtN map. To define this map for the Helmholtz problem we first require the Helmholtz extension operator from the boundary of a subdomain $\Omega_{s}$.

Let $\Gamma_{s} = \partial\Omega_{s} \setminus \partial\Omega$ and suppose we have Dirichlet data $v_{\Gamma_{s}}$ on $\Gamma_{s}$, then the Helmholtz extension $v$ in $\Omega_{s}$ is defined as the solution of
\begin{subequations}
\label{HelmholtzExtension}
\begin{align}
\label{HelmholtzExtensionEquation}
-\Delta v - k^{2} v & = 0 & & \text{in } \Omega_{s},\\
\label{HelmholtzExtensionDirichletBC}
v & = v_{\Gamma_{s}} & & \text{on } \Gamma_{s},\\
\label{HelmholtzExtensionProblemBC}
\mathcal{C}(v) & = 0 & & \text{on } \partial\Omega_{s}\cap\partial\Omega,
\end{align}
\end{subequations}
where $\mathcal{C}(v) = 0$ represents the original problem boundary conditions, as in \eqref{HelmholtzBCs}. The DtN map takes in the Dirichlet data $v_{\Gamma_{s}}$ on $\Gamma_{s}$ and gives as output the corresponding Neumann data, that is
\begin{align}
\label{DtNMap}
\mathrm{DtN}_{\Omega_{s}}(v_{\Gamma_{s}}) = \left.\frac{\partial v}{\partial \boldsymbol{n}_s} \right\rvert_{\Gamma_{s}}
\end{align}
where $v$ is the Helmholtz extension defined by \eqref{HelmholtzExtension}.

We now seek eigenfunctions of the DtN map locally on each subdomain $\Omega_{s}$, given by solving
\begin{align}
\label{DtNEigenproblem}
\mathrm{DtN}_{\Omega_{s}}(u_{\Gamma_{s}}) = \lambda u_{\Gamma_{s}}
\end{align}
for eigenfunctions $u_{\Gamma_{s}}$ and eigenvalues $\lambda \in \mathbb{C}$. To provide functions to go into the coarse space, we take the Helmholtz extension of $u_{\Gamma_{s}}$ in $\Omega_{s}$ and then extend by zero into the whole domain $\Omega$ using the partition of unity.

To formulate the discrete problem, we require the coefficient matrices $\widetilde{A}_{s}$ corresponding to local Neumann problems on $\Omega_{s}$ with boundary conditions $\mathcal{C}=0$ on $\partial\Omega_{s}\cap\partial\Omega$, defined analogously to that of the local Robin problems in \eqref{ORASLocalSystem}. Further, we need to distinguish between degrees of freedom on the boundary $\Gamma_{s}$ and the interior of the subdomain $\Omega_{s}$, as such we let $\Gamma_{s}$ and $I_{s}$ be the set of indices on the boundary and interior respectively. We also let
\begin{align}
\label{DtNMassMatrix}
M_{\Gamma_{s}} = \left(\int_{\Gamma_{s}} \phi_{j} \phi_{i} \right)_{i,j \in \Gamma_{s}}
\end{align}
denote the mass matrix on the subdomain interface. Using standard block notation to denote submatrices of $A_{s}$ and $\widetilde{A}_{s}$ the discrete DtN eigenproblem is
\begin{align}
\label{DiscreteDtNEigenproblem}
\left(\widetilde{A}_{\Gamma_{s},\Gamma_{s}} - A_{\Gamma_{s},I_{s}}A_{I_{s},I_{s}}^{-1}A_{I_{s},\Gamma_{s}}\right) \mathbf{u}_{\Gamma_{s}} = \lambda M_{\Gamma_{s}} \mathbf{u}_{\Gamma_{s}}.
\end{align}
The Helmholtz extension of $\mathbf{u}_{\Gamma_{s}}$ to degrees of freedom in $I_{s}$ is then given by $\mathbf{u}_{I_{s}} = - A_{I_{s},I_{s}}^{-1}A_{I_{s},\Gamma_{s}} \mathbf{u}_{\Gamma_{s}}$. Letting $\mathbf{u}_{s}$ denote the Helmholtz extension, the corresponding vector which enters the coarse space $Z$ is $R_{s}^{T} D_{s} \mathbf{u}_{s}$. For further details and motivation behind the DtN eigenproblems see \cite{Conen:2014:ACS}.

It remains to determine which eigenfunctions of \eqref{DiscreteDtNEigenproblem} should go on to be included in the coarse space. Several selection criteria were investigated in \cite{Conen:2014:ACS} and it was clear that the best choice was to select eigenvectors corresponding to eigenvalues with the smallest real part. That is, we use a threshold on the abscissa $\eta = \mathrm{Re}(\lambda)$ given by
\begin{align}
\label{Threshold}
\eta < \eta_{\text{max}},
\end{align}
where $\eta_{\text{max}}$ depends on $k_{s} = \max_{\vec{x}\in\Omega_{s}} k(\vec{x})$. In particular, the choice $\eta_{\text{max}} = k_{s}$ is advocated in \cite{Conen:2014:ACS}. Note that for larger $\eta_{\text{max}}$ more eigenfunctions are included in the coarse space, increasing its size and the associated computational cost. Nonetheless, it was recently showed that taking a slightly larger threshold $\eta_{\text{max}} = k_{s}^{4/3}$ can be beneficial in certain cases in order to gain robustness to the wave number \cite{Bootland:2019:ODN}. However, this only occurs for the homogeneous problem with sufficiently uniform subdomains. Since it is not necessarily known in advance how many eigenvalues are below the threshold and in order to make a fair comparison in our numerical tests, we will consider using a fixed number of eigenvectors per subdomain.

\subsection{The GenEO coarse space}
\label{subsec:GenEO}

The GenEO (Generalised Eigenproblems in the Overlap) coarse space was derived in \cite{Spillane:2014:ARC} to provide a rigorously robust approach for symmetric positive definite problems even in the presence of heterogeneities. The fundamental generalised eigenproblems on $\Omega_{s}$ at the variational level are given by
\begin{align}
\label{GenEOVariationalEigenproblem}
a_{\Omega_{s}}(u,v) & = \lambda a_{\Omega_{s}^{\circ}}\left(\Xi_{s}(u),\Xi_{s}(v)\right) & & \forall \ v \in V(\Omega_{s}),
\end{align}
where $\Xi_{s}$ represents the action of the partition of unity on $\Omega_{s}$ and $\Omega_{s}^{\circ}$ is the overlapping zone, that is the part of $\Omega_{s}$ which overlaps with any other subdomain. Here $a_{D}(\cdot,\cdot)$ stems from the underlying variational problem on the domain $D$, in particular with problem boundary conditions on $\partial\Omega$ and natural (Neumann) conditions on parts of $\partial D$ internal to $\Omega$. The particular form of eigenproblem in \eqref{GenEOVariationalEigenproblem} arises naturally in the analysis of \cite{Spillane:2014:ARC}. We also note that \eqref{GenEOVariationalEigenproblem} possesses infinite eigenvalues when there exists $u \neq 0$ such that $a_{\Omega_{s}^{\circ}}\left(\Xi_{s}(u),\Xi_{s}(v)\right) = 0$  $\forall \ v \in V(\Omega_{s})$ but $a_{\Omega_{s}}(u,v) \neq 0$ for some $v \in V(\Omega_{s})$, for example when $u$ is supported only outside $\Omega_{s}^{\circ}$.

The discrete form of the eigenproblem \eqref{GenEOVariationalEigenproblem} is
\begin{align}
\label{DiscreteGenEOEigenproblem}
\widetilde{A}_{s} \mathbf{u} = \lambda D_{s} \widetilde{A}_{s}^{\circ} D_{s} \mathbf{u}
\end{align}
where $\widetilde{A}_{s}^{\circ}$ is the (Neumann) matrix built from assembling only over elements in the overlapping zone $\Omega_{s}^{\circ}$. The eigenfunctions selected to enter the coarse space are the low frequency modes, that is those corresponding to the smallest eigenvalues. Typically, either a fixed number of eigenfunctions are taken per subdomain or a threshold $\lambda < \lambda_{\text{max}}$ on the corresponding eigenvalues is used, where $\lambda_{\text{max}}$ can be chosen based on problem parameters to achieve a specified condition number of the preconditioned system; see \cite{Spillane:2014:ARC}.

The precise formulation of the eigenproblem and use of the overlapping zone is somewhat flexible. In particular, the overlapping zone can be replaced with the whole subdomain, as in \cite{Haferssas:2017:ADS} which utilises an eigenproblem of the form
\begin{align}
\label{DiscreteGenEO2Eigenproblem}
\widetilde{A}_{s} \mathbf{u} = \lambda D_{s} A_{s} D_{s} \mathbf{u}.
\end{align}
It is this form of GenEO that we shall build upon shortly. We note that two separate GenEO eigenproblems can also be formulated to provide bounds on both ends of the preconditioned spectrum, as is found when using a symmetric ORAS approach in \cite{Haferssas:2017:ADS}. This flexibility and robustness of GenEO-type methods has yet to be fully explored, especially for problems which are not symmetric positive definite where current theory breaks down. We now consider the utility of using GenEO as basis for constructing a coarse space for the heterogeneous Helmholtz problem.

\subsection{A GenEO-type coarse space for the Helmholtz equation}
\label{subsec:GenEOForHelmholtz}

In pursuing a GenEO approach for the Helmholtz equation we must first note the matrices $\widetilde{A}_{s}$ and $A_{s}$, now stemming from the bilinear form \eqref{WeakFormBilinearPart}, are no longer symmetric positive definite. As such, eigenvalues $\lambda$ of the eigenproblem \eqref{DiscreteGenEO2Eigenproblem} are no longer real and positive (or infinite). As a threshold criterion, as with the DtN approach, we can consider the abscissa $\eta = \mathrm{Re}(\lambda)$ instead and seek eigenfunctions corresponding to $\eta < \eta_{\text{max}}$. Further, since the theory breaks down without the symmetric positive definite assumption, it is no longer clear that \eqref{DiscreteGenEO2Eigenproblem} provides appropriate eigenvectors for the coarse space. Indeed, applying out of the box the GenEO method using \eqref{DiscreteGenEO2Eigenproblem} fails to provide a satisfactory method.

Since GenEO is designed for positive definite problems, a natural proposal is to use a nearby positive definite problem in the formulation of the eigenproblem, namely a Laplace problem, that is setting $k=0$ for the purpose of constructing the coarse space. Let $L_{s}$ be local Dirichlet matrix corresponding to the discrete Laplacian and $\widetilde{L}_{s}$ the equivalent Neumann matrix on $\Omega_{s}$, then
\begin{align}
\label{DiscreteLaplaceGenEOEigenproblem}
\widetilde{L}_{s} \mathbf{u} = \lambda D_{s} L_{s} D_{s} \mathbf{u}.
\end{align}
has positive real eigenvalues. While this approach can perform reasonably well when the wave number $k$ is small, the behaviour as $k$ grows becomes increasingly poor, as might be expected given the coarse space is independent of $k$.

To incorporate $k$, we instead link the underlying Helmholtz problem to the positive definite Laplace problem to formulate a GenEO-type eigenproblem
\begin{align}
\label{DiscreteH-GenEOEigenproblem}
\widetilde{A}_{s} \mathbf{u} = \lambda D_{s} L_{s} D_{s} \mathbf{u}.
\end{align}
Eigenvalues of \eqref{DiscreteH-GenEOEigenproblem} are now, in general, complex (though we note they primarily appear to cluster close to the real line) and so we suggest to threshold based on the abscissa $\eta < \eta_{\text{max}}$. We call this GenEO-type approach for the Helmholtz problem ``H-GenEO''. Some initial exploration of this method can be found in \cite{Bootland:HGenEOPresentation}, where the approach is seen to perform well for a 2D wave-guide problem and provide robustness to heterogeneity as well as to the wave number $k$, albeit requiring a comparatively large coarse space. Again, in our numerical experiments we take a fixed number of eigenvectors and vary this quantity, aiming to give a relatively fair comparison of the spectral coarse space methods.

\section{Implementation details}
\label{sec:ImplementationDetails}

Numerical results within this paper have been obtained using
FreeFEM~\cite{hecht2012new}. More precisely, ffddm---a light-weight layer in
FreeFEM domain specific language that generates all domain decomposition data
structures---was used on top of HPDDM~\cite{jolivet2013scalable}, which handles
the underlying computations such as matrix--vector products or preconditioner
applications. When it comes to spectral two-level methods, as depicted in the
previous \Cref{subsec:DtN,subsec:GenEO,subsec:GenEOForHelmholtz}, the bulk of
the work comes from the local eigensolves, see, e.g.,
\cref{DiscreteDtNEigenproblem,DiscreteGenEOEigenproblem,DiscreteH-GenEOEigenproblem}.
These are performed concurrently in each subdomain, using solvers such as
ARPACK~\cite{lehoucq1998arpack} or SLEPc~\cite{hernandez2005ssf}. Then, the construction of the algebraic Galerkin operator $Q$
introduced in \Cref{sec:DomainDecomposition} is performed in HPDDM. In order to
deal efficiently with both large numbers of subdomains and large numbers of
local eigenvectors, the assembly and exact $LU$ factorisation of $Q$ is
performed on a subset of processes from the global communicator. We use MUMPS~\cite{amestoy2001fully}
both as a subdomain solver and a coarse operator solver in our experiments.
Again, we refer interested readers to~\cite{jolivet2013scalable} for more
details about this subject.

\begin{figure}
\centering
    \subfloat[Coarse mesh\label{fig:coarse-mesh}]{\includegraphics[width=0.4\textwidth]{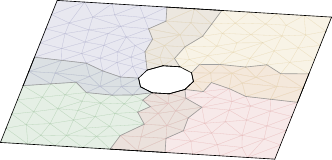}}\\
    \subfloat[Refined mesh with minimal overlap\label{fig:fine-mesh-min}]{\includegraphics[width=0.4\textwidth]{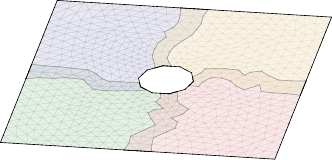}}\qquad
    \subfloat[Refined mesh with coarse overlap\label{fig:fine-mesh-coarse}]{\includegraphics[width=0.4\textwidth]{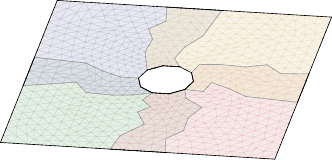}}
    \caption{Different overlap definitions (bottom) for a uniformly refined simple 2D coarse mesh (top).\label{fig:coarse-fine}}
\end{figure}

For the grid coarse space method from \Cref{subsec:CoarseGrid}, two meshes of
the same domain must be considered simultaneously. First a coarse one, see an
example~\Cref{fig:coarse-mesh}, and then a fine one, see
examples~\Cref{fig:fine-mesh-min,fig:fine-mesh-coarse}. We generate fine meshes
using a uniform refinement in which edges of all triangles or tetrahedra are
divided uniformly in $s$. In the aforementioned example, $s=2$. Note that if
one wants $n_\text{ppwl} h \approx \lambda$ on the fine grid, on which the
original algebraic system of equations~\eqref{LinearSystem} stems from, then the
coarse mesh has to satisfy $H_{\text{coarse}} = \frac{h}{s}$. When it comes to
the underlying domain decomposition preconditioners, two options will be
considered next. One where both levels of overlap on fine and coarse grids are
minimal, i.e.~1, see~\Cref{fig:fine-mesh-min}, and another where the level of
overlap on the coarse grid is minimal, i.e.~1, and where the level of overlap
on the fine grid is $s$, see~\Cref{fig:fine-mesh-coarse}. In the latter case,
note that subdomains on the fine grid are merely uniformly refined subdomains
of the coarse grid. While the setup cost of such a two-level preconditioner is
much lower than for a spectral method as described in the previous paragraph,
applying it to a vector, for example in a Krylov method, is usually costlier.
Indeed, the coarse problem is in this case often solved iteratively, instead of
using an exact factorisation of the coarse operator. This outer--inner strategy,
which also makes the use of flexible methods such as FGMRES~\cite{saad1993flexible} mandatory,
may not perform very well prior to some tuning of the inner (coarse) solver.

In the numerical experiments presented in this paper, the inner coarse problem in the grid coarse space method is defined with a splitting level $s=2$ and is solved approximately with GMRES preconditioned by a one-level method, with a tolerance of $0.1$. Additionally, in the spirit of ``shifted Laplacian preconditioning'' (see~\cite{Erlangga:2004:OCP,lahaye2017modern}), the coarse problem is defined with added absorption in the equation. This improves the convergence of the one-level method, which then requires fewer iterations to reach the prescribed inner tolerance. Two-level domain decomposition preconditioners with added absorption have been proposed in the literature for the Helmholtz~\cite{Graham:2017:DDP} and Maxwell~\cite{Bonazzoli:2019:ADD} equations. The amount of added absorption needs to be chosen carefully: if the amount of added absorption is too small, there is no gain in the convergence of the inner coarse problem. On the other hand, if the amount is too large, the coarse problem becomes a bad approximation of the original problem (without absorption), and the number of outer iterations increases. We choose the additional imaginary term to be proportional to the wave number $k$ as a compromise.

We conclude this section by mentioning that all the spectral preconditioners may be used in PETSc~\cite{petsc-user-ref} through the PCHPDDM preconditioner~\cite{jolivet2020petsc}. When it comes to the grid coarse space method, we have a custom implementation in HPDDM that handles systems with multiple right-hand sides, which are encountered frequently in wave propagation phenomena~\cite{tournier2016micro,roux2017block,Dolean:2020:IFD}. PCMG, the PETSc machinery for geometric multigrid, does not handle such systems as of version 3.14.1.

\section{Comparative numerical studies}
\label{sec:ComparativeNumericalStudies}

In this section we compare the domain decomposition preconditioners for several challenging heterogeneous test problems. We include problems in both 2D and 3D and vary different characteristics of the problem to understand how they can affect performance. In particular, we are interested in the behaviour for larger wave numbers (or equivalently larger frequencies) and for an increasing number of subdomains within the domain decomposition. We will also consider discretisations with differing number of points per wavelength and will see how this can affect performance of the preconditioners. For the approaches that use a spectral coarse space, we further investigate the choice of how many eigenvectors should be taken per subdomain.

In the results that follow, a tabulated value of $-$ signifies that the particular test instance did not run, typically this is for the smallest problem when trying to solve on too many subdomains. When the maximum number of iterations is reached before convergence, we denote this by $\times$. Depending on the test problem, the maximum number of iterations is set at either 500 or 1000.

\subsection{The Marmousi problem}
\label{sec:MarmousiProblem}

\begin{figure}[H]
	\centering
	\includegraphics[width=\textwidth,trim=0cm 0cm 2cm 0cm ,clip]{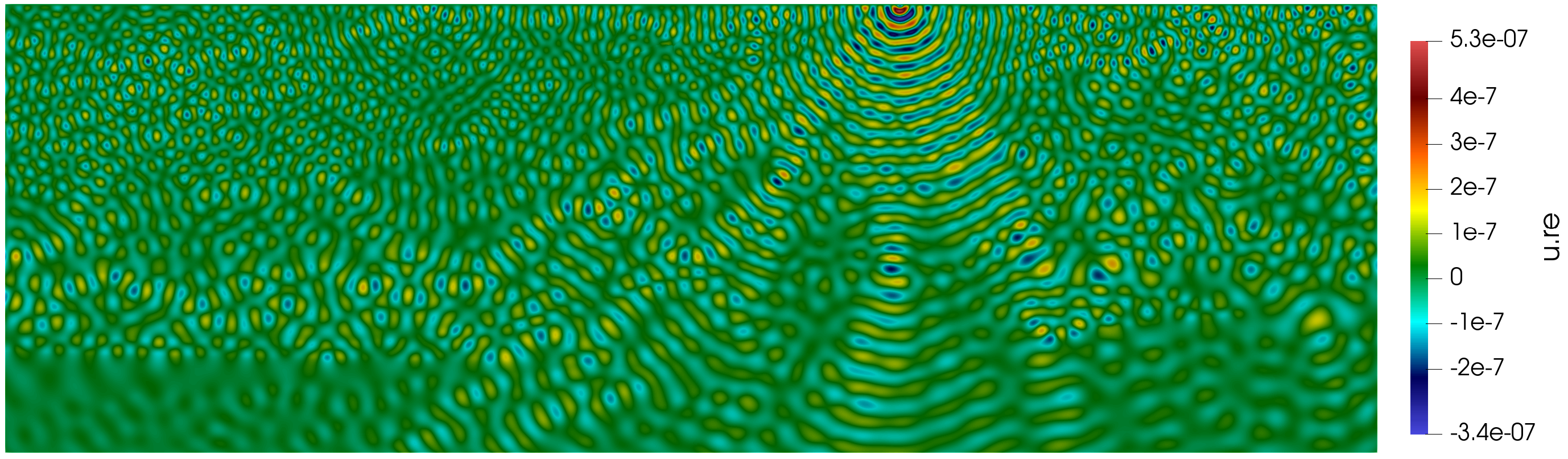}
	\caption{The real part of the solution to the Marmousi problem at a frequency of $20\,$Hz.}
	\label{fig:marmousi}
\end{figure}

The Marmousi problem, see \cite{versteeg1990marmousi}, has become a benchmark geophysics test case for the assessment of numerical methods and solvers when used in a direct or inverse problem context. It features a 2D rectangular domain with velocity data modelling a heterogeneous subsurface. We utilise this as the wave speed within our Helmholtz model with a point source placed close to the surface; see the example solution in \Cref{fig:marmousi}. The problem is posed on the domain $(0,9.2)\times(-3,0)$ with homogeneous Dirichlet conditions imposed on the top (surface) boundary and Robin conditions on the remaining (subsurface) boundaries, the latter being artificial boundaries set for computational purposes.

The difficulty of this test case stems both from the heterogeneous nature of the problem but also from the presence of a substantial number of wavelengths and for this reason, the problem is considered to be high-frequency. Therefore, it is important to have a sufficiently accurate discretisation, keeping in mind that for inverse problems the accuracy requirements are different (usually 5 points per wavelength) and weaker than in the case of only solving the forward problem (in common engineering practice, we use 10 points per wavelength).

We first consider discretisation with 5 points per wavelength. In \Cref{Table:Marmousi-P2-nppwl5-non-spectral} we show results for the one-level method with minimal overlap and also with minimal overlap from the coarse grid (twice the overlap---see \Cref{fig:coarse-fine}) as references to compare to. Further, \Cref{Table:Marmousi-P2-nppwl5-non-spectral} gives results for the two-level coarse grid approach with minimal overlap from the coarse grid. We see that both cases of the one-level method are not scalable (iteration counts increase as the number of subdomains $N$ increases) and perform worse as the frequency $f$ becomes higher. On the other hand, the two-level coarse grid approach is reasonably scalable and with only mildly increasing iteration counts as the frequency increases.

\begin{table}[H]
	\centering
	\caption{Results using the one-level and coarse grid methods for the Marmousi problem when using 5 points per wavelength, varying the frequency $f$ and the number of subdomains $N$.}
	\label{Table:Marmousi-P2-nppwl5-non-spectral}
	\tabulinesep=1.2mm
	\footnotesize
	\begin{tabu}{c|ccccc|ccccc|ccccc}
		& \multicolumn{5}{c|}{One-level (min.\ overlap)} & \multicolumn{5}{c|}{One-level (coarse overlap)} & \multicolumn{5}{c}{Coarse grid} \\
		\hline
		$f \ \backslash \ N$ & 10 & 20 & 40 & 80 & 160 & 10 & 20 & 40 & 80 & 160 & 10 & 20 & 40 & 80 & 160 \\
		\hline
		1 & 30 & 46 & 78 & 98 & $-$ & 26 & 39 & 47 & 64 & $-$ & 15 & 18 & 19 & 20 & $-$ \\
		5 & 57 & 82 & 117 & 170 & 236 & 53 & 76 & 105 & 154 & 213 & 26 & 29 & 28 & 29 & 31 \\
		10 & 75 & 111 & 173 & 232 & 330 & 68 & 102 & 158 & 212 & 302 & 32 & 35 & 41 & 40 & 42 \\
		20 & 89 & 133 & 193 & 268 & 373 & 82 & 125 & 178 & 248 & 347 & 34 & 35 & 42 & 43 & 44
	\end{tabu}
\end{table}

We now turn to the two spectral coarse space methods with results in \Cref{Table:Marmousi-P2-nppwl5-spectral-nu} being given for a varying number of eigenvectors taken per subdomain $\nu$. We first see that the H-GenEO method within this scenario performs poorly and, without sufficiently many eigenvectors being taken, often exhibits iteration counts worse than the one-level method. On the other hand, the DtN method generally converges faster than the one-level method and, with enough eigenvectors, can give reasonable performance. However, as the frequency increases the approach begins to struggle a little with larger iteration counts typically required for convergence even with more eigenvectors being utilised. Nonetheless, for small frequencies the approach can beat the coarse grid method with relatively few eigenvectors per subdomain and even for the larger frequencies tested, when using the most subdomains and sufficient eigenvectors per subdomain the DtN approach converges faster than the coarse grid method. We note by way of caution though that taking too many eigenvectors per subdomain can eventually lead to some deterioration in the performance of the DtN method. We can also see from \Cref{Table:Marmousi-P2-nppwl5-spectral-nu} that the larger coarse overlap is beneficial for the DtN method and primarily yields smaller iterations counts. On the other hand, minimal overlap appears preferable for H-GenEO.

\begin{sidewaystable}[ht]
	\centering
	\caption{Results using the H-GenEO and DtN spectral coarse space methods for the Marmousi problem when using 5 points per wavelength, varying the frequency $f$, the number of subdomains $N$, and the number of eigenvectors used per subdomain $\nu$.}
	\label{Table:Marmousi-P2-nppwl5-spectral-nu}
	\tabulinesep=1.2mm
	\footnotesize
	\begin{tabu}{cc|ccccc|ccccc|ccccc|ccccc}
		& & \multicolumn{5}{c|}{H-GenEO (min.\ overlap)} & \multicolumn{5}{c|}{DtN (min.\ overlap)} & \multicolumn{5}{c|}{H-GenEO (coarse overlap)} & \multicolumn{5}{c}{DtN (coarse overlap)} \\
		\hline
		$\nu$ & $f \ \backslash \ N$ & 10 & 20 & 40 & 80 & 160 & 10 & 20 & 40 & 80 & 160 & 10 & 20 & 40 & 80 & 160 & 10 & 20 & 40 & 80 & 160 \\
		\hline
		\multirow{4}{*}{20} & 1 & 9 & 14 & 21 & 21 & $-$ & 8 & 7 & 6 & 6 & $-$ & 16 & 22 & 31 & 34 & $-$ & 5 & 5 & 5 & 5 & $-$ \\
		& 5 & 46 & 53 & 67 & 76 & 100 & 37 & 46 & 59 & 74 & 85 & 58 & 86 & 120 & 158 & 201 & 39 & 54 & 74 & 83 & 68 \\
		& 10 & 67 & 99 & 146 & 173 & 209 & 61 & 87 & 118 & 156 & 190 & 76 & 133 & 220 & 301 & 420 & 62 & 91 & 129 & 178 & 242 \\
		& 20 & 81 & 120 & 169 & 248 & 378 & 75 & 121 & 157 & 208 & 275 & 83 & 125 & 226 & 370 & $\times$ & 75 & 117 & 158 & 216 & 300 \\
		\hline
		\multirow{4}{*}{40} & 1 & 9 & 13 & 16 & 21 & $-$ & 10 & 8 & 6 & 6 & $-$ & 16 & 22 & 26 & 32 & $-$ & 7 & 8 & 6 & 7 & $-$ \\
		& 5 & 31 & 32 & 42 & 49 & 71 & 25 & 34 & 36 & 40 & 37 & 46 & 74 & 104 & 137 & 186 & 28 & 32 & 24 & 11 & 9 \\
		& 10 & 66 & 85 & 105 & 106 & 126 & 48 & 65 & 90 & 102 & 122 & 100 & 153 & 233 & 269 & 373 & 52 & 80 & 115 & 131 & 106 \\
		& 20 & 86 & 136 & 209 & 264 & 323 & 68 & 106 & 138 & 186 & 259 & 119 & 232 & 391 & $\times$ & $\times$ & 68 & 105 & 149 & 236 & 357 \\
		\hline
		\multirow{4}{*}{80} & 1 & 9 & 12 & 18 & 21 & $-$ & 5 & 6 & 6 & 6 & $-$ & 15 & 20 & 22 & 33 & $-$ & 4 & 4 & 6 & 6 & $-$ \\
		& 5 & 20 & 20 & 28 & 41 & 62 & 15 & 18 & 13 & 39 & 67 & 40 & 62 & 87 & 127 & 184 & 10 & 8 & 8 & 23 & 19 \\
		& 10 & 47 & 56 & 64 & 72 & 88 & 33 & 40 & 55 & 48 & 39 & 88 & 118 & 171 & 225 & 327 & 39 & 46 & 42 & 14 & 13 \\
		& 20 & 87 & 125 & 171 & 186 & 217 & 54 & 86 & 111 & 152 & 164 & 153 & 246 & 404 & $\times$ & $\times$ & 57 & 104 & 155 & 188 & 145 \\
		\hline
		\multirow{4}{*}{160} & 1 & 9 & 11 & 17 & 21 & $-$ & 8 & 7 & 5 & 6 & $-$ & 16 & 15 & 21 & 33 & $-$ & 6 & 7 & 9 & 6 & $-$ \\
		& 5 & 15 & 17 & 26 & 37 & 56 & 13 & 25 & 23 & 26 & 36 & 37 & 59 & 82 & 116 & 171 & 7 & 19 & 10 & 8 & 19 \\
		& 10 & 33 & 40 & 45 & 56 & 73 & 21 & 27 & 30 & 48 & 80 & 68 & 94 & 143 & 196 & 284 & 18 & 19 & 17 & 48 & 29 \\
		& 20 & 64 & 83 & 121 & 134 & 157 & 33 & 59 & 67 & 76 & 75 & 131 & 192 & 320 & 403 & $\times$ & 43 & 75 & 77 & 61 & 35
	\end{tabu}
\end{sidewaystable}

Overall, while the DtN approach can reduce the number of iterations required for convergence compared to the coarse grid method in some cases, particularly low frequencies, the coarse grid approach exhibits greater robustness and requires relatively low iterations counts, suggesting it is the favourable approach in this scenario of using 5 points per wavelength.

We now consider the case of discretisation with 10 points per wavelength and similarly compare the different approaches. In \Cref{Table:Marmousi-P2-nppwl10-non-spectral} we show results for the one-level methods along with the coarse grid method. We notice that the iteration counts in this situation are slightly higher than when using 5 points per wavelength but otherwise the picture remains similar with the coarse grid approach giving a reasonably robust method, albeit with some slow growth in iteration counts as the frequency is increased.

\begin{table}[H]
	\centering
	\caption{Results using the one-level and coarse grid methods for the Marmousi problem when using 10 points per wavelength, varying the frequency $f$ and the number of subdomains $N$.}
	\label{Table:Marmousi-P2-nppwl10-non-spectral}
	\tabulinesep=1.2mm
	\footnotesize
	\begin{tabu}{c|ccccc|ccccc|ccccc}
		& \multicolumn{5}{c|}{One-level (min.\ overlap)} & \multicolumn{5}{c|}{One-level (coarse overlap)} & \multicolumn{5}{c}{Coarse grid} \\
		\hline
		$f \ \backslash \ N$ & 10 & 20 & 40 & 80 & 160 & 10 & 20 & 40 & 80 & 160 & 10 & 20 & 40 & 80 & 160 \\
		\hline
		1 & 34 & 49 & 72 & 143 & $-$ & 30 & 43 & 63 & 97 & $-$ & 16 & 18 & 19 & 21 & $-$ \\
		5 & 62 & 94 & 137 & 191 & 268 & 58 & 87 & 126 & 175 & 246 & 29 & 29 & 34 & 34 & 36 \\
		10 & 85 & 136 & 185 & 272 & 371 & 78 & 124 & 172 & 251 & 346 & 35 & 41 & 43 & 46 & 45 \\
		20 & 101 & 152 & 213 & 299 & 419 & 92 & 142 & 198 & 272 & 389 & 39 & 47 & 48 & 49 & 49
	\end{tabu}
\end{table}

Turning to the spectral coarse space approaches, we see in \Cref{Table:Marmousi-P2-nppwl10-spectral-nu} that H-GenEO now becomes competitive. For larger frequencies and equivalent numbers of eigenvectors taken per subdomain, the H-GenEO method generally requires fewer iterations and this is particularly true for the large frequencies with many subdomains. We note that the minimal overlap case is preferable for H-GenEO while for DtN there is no longer a clear-cut preference between the coarse overlap and minimal overlap cases. Provided sufficiently many eigenvectors are taken per subdomain, we see that H-GenEO now outperforms the coarse grid approach in terms of iterations counts and is particularly strong for the large frequency and many subdomains situation.

Overall, in the scenario of 10 points per wavelength, we have different findings from the 5 point per wavelength scenario. Now, while the coarse grid approach is still relatively robust, the H-GenEO method can be seen as strong alternative which is able to markedly reduce the number of iterations required in the case of large frequencies and many subdomains. One possible explanation is that spectral information is not relevant or sufficiently well approximated in the case of a low number of points per wavelength.

\begin{sidewaystable}[ht]
	\centering
	\caption{Results using the H-GenEO and DtN spectral coarse space methods for the Marmousi problem when using 10 points per wavelength, varying the frequency $f$, the number of subdomains $N$, and the number of eigenvectors used per subdomain $\nu$.}
	\label{Table:Marmousi-P2-nppwl10-spectral-nu}
	\tabulinesep=1.2mm
	\footnotesize
	\begin{tabu}{cc|ccccc|ccccc|ccccc|ccccc}
		& & \multicolumn{5}{c|}{H-GenEO (min.\ overlap)} & \multicolumn{5}{c|}{DtN (min.\ overlap)} & \multicolumn{5}{c|}{H-GenEO (coarse overlap)} & \multicolumn{5}{c}{DtN (coarse overlap)} \\
		\hline
		$\nu$ & $f \ \backslash \ N$ & 10 & 20 & 40 & 80 & 160 & 10 & 20 & 40 & 80 & 160 & 10 & 20 & 40 & 80 & 160 & 10 & 20 & 40 & 80 & 160 \\
		\hline
		\multirow{4}{*}{20} & 1 & 8 & 8 & 9 & 15 & $-$ & 8 & 8 & 5 & 7 & $-$ & 7 & 7 & 9 & 14 & $-$ & 5 & 6 & 4 & 4 & $-$ \\
		& 5 & 42 & 42 & 39 & 32 & 19 & 42 & 49 & 59 & 51 & 56 & 40 & 46 & 47 & 38 & 32 & 39 & 49 & 67 & 72 & 85 \\
		& 10 & 78 & 112 & 125 & 133 & 99 & 70 & 105 & 128 & 157 & 166 & 72 & 108 & 133 & 156 & 138 & 65 & 99 & 128 & 169 & 197 \\
		& 20 & 96 & 136 & 180 & 259 & 342 & 92 & 131 & 175 & 232 & 318 & 89 & 130 & 177 & 248 & 345 & 85 & 125 & 170 & 227 & 315 \\
		\hline
		\multirow{4}{*}{40} & 1 & 7 & 8 & 9 & 14 & $-$ & 6 & 5 & 6 & 6 & $-$ & 7 & 7 & 9 & 13 & $-$ & 4 & 3 & 3 & 4 & $-$ \\
		& 5 & 22 & 19 & 17 & 13 & 14 & 25 & 29 & 31 & 32 & 33 & 21 & 20 & 19 & 17 & 20 & 25 & 34 & 43 & 40 & 38 \\
		& 10 & 58 & 63 & 54 & 42 & 26 & 49 & 71 & 74 & 84 & 82 & 58 & 70 & 65 & 61 & 44 & 48 & 75 & 88 & 112 & 128 \\
		& 20 & 101 & 148 & 169 & 161 & 148 & 85 & 120 & 159 & 190 & 217 & 94 & 141 & 183 & 203 & 214 & 80 & 116 & 154 & 194 & 256 \\
		\hline
		\multirow{4}{*}{80} & 1 & 7 & 8 & 9 & 13 & $-$ & 7 & 7 & 5 & 5 & $-$ & 6 & 7 & 9 & 12 & $-$ & 5 & 4 & 3 & 4 & $-$ \\
		& 5 & 13 & 12 & 11 & 10 & 12 & 16 & 16 & 19 & 17 & 10 & 12 & 11 & 12 & 13 & 17 & 19 & 15 & 14 & 7 & 6 \\
		& 10 & 31 & 28 & 23 & 17 & 15 & 28 & 39 & 40 & 42 & 41 & 30 & 31 & 27 & 23 & 24 & 29 & 47 & 55 & 66 & 34 \\
		& 20 & 81 & 90 & 78 & 53 & 39 & 63 & 89 & 93 & 100 & 107 & 85 & 104 & 107 & 89 & 69 & 64 & 93 & 108 & 147 & 202 \\
		\hline
		\multirow{4}{*}{160} & 1 & 7 & 8 & 8 & 13 & $-$ & 4 & 7 & 5 & 6 & $-$ & 6 & 7 & 8 & 11 & $-$ & 3 & 6 & 3 & 5 & $-$ \\
		& 5 & 10 & 9 & 10 & 10 & 12 & 10 & 11 & 12 & 17 & 24 & 8 & 9 & 10 & 13 & 16 & 10 & 9 & 8 & 21 & 14 \\
		& 10 & 20 & 16 & 14 & 13 & 13 & 19 & 23 & 25 & 25 & 24 & 18 & 15 & 16 & 16 & 19 & 21 & 29 & 30 & 25 & 17 \\
		& 20 & 45 & 40 & 34 & 25 & 19 & 35 & 46 & 47 & 56 & 59 & 50 & 52 & 50 & 39 & 28 & 36 & 61 & 77 & 84 & 76
	\end{tabu}
\end{sidewaystable}

\subsection{The COBRA cavity}
\label{sec:CobraCavity}

\begin{figure}[H]
	\centering
	\subfloat[3D view\label{fig:cobra-3d-view}]{\includegraphics[width=0.48\textwidth,trim=2cm 0cm 2cm 0cm ,clip]{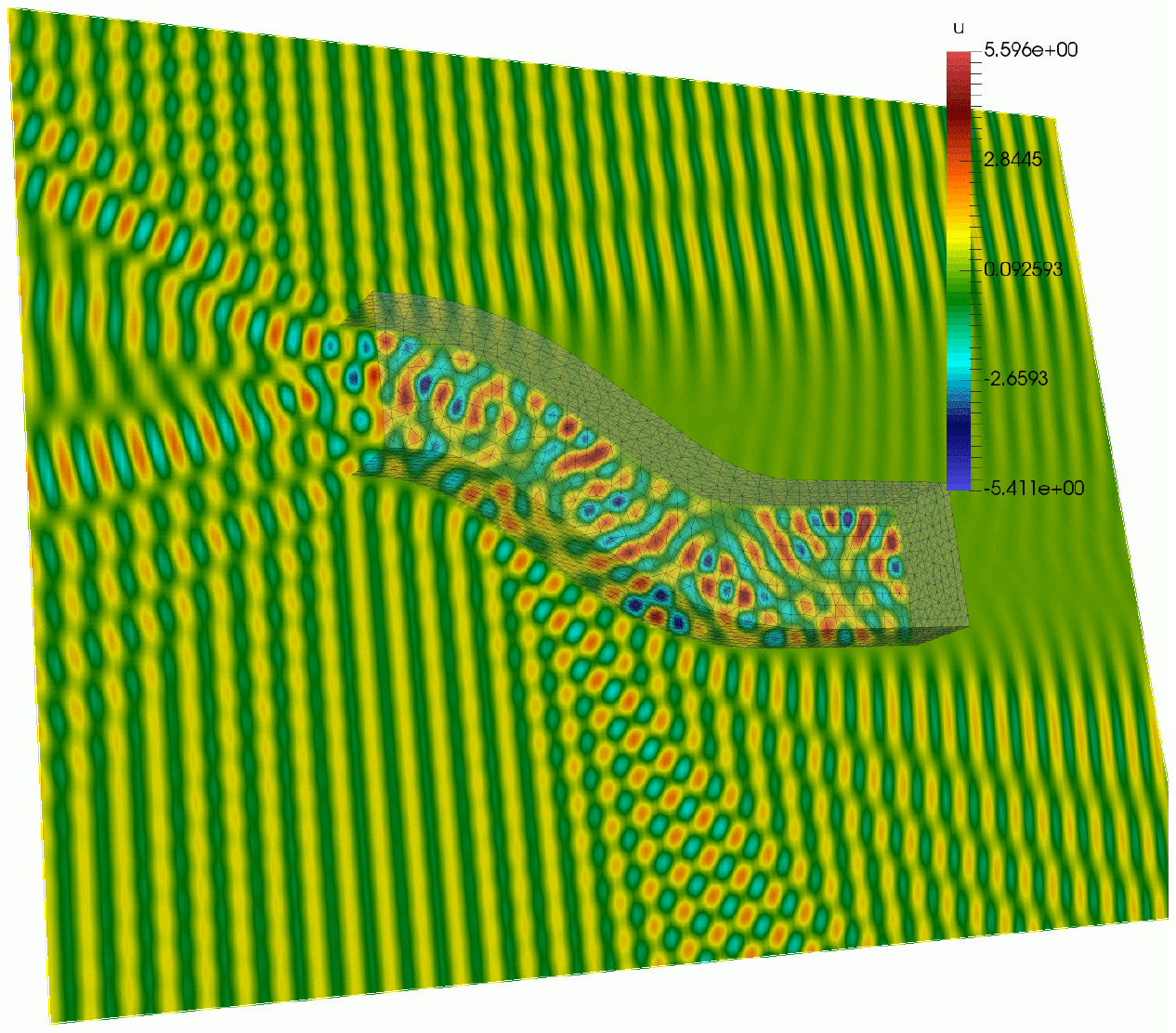}}
	\hfill
	\subfloat[2D cross-section\label{fig:cobra-2d-cross-section}]{\includegraphics[width=0.48\textwidth,trim=2cm 0cm 2cm 0cm ,clip]{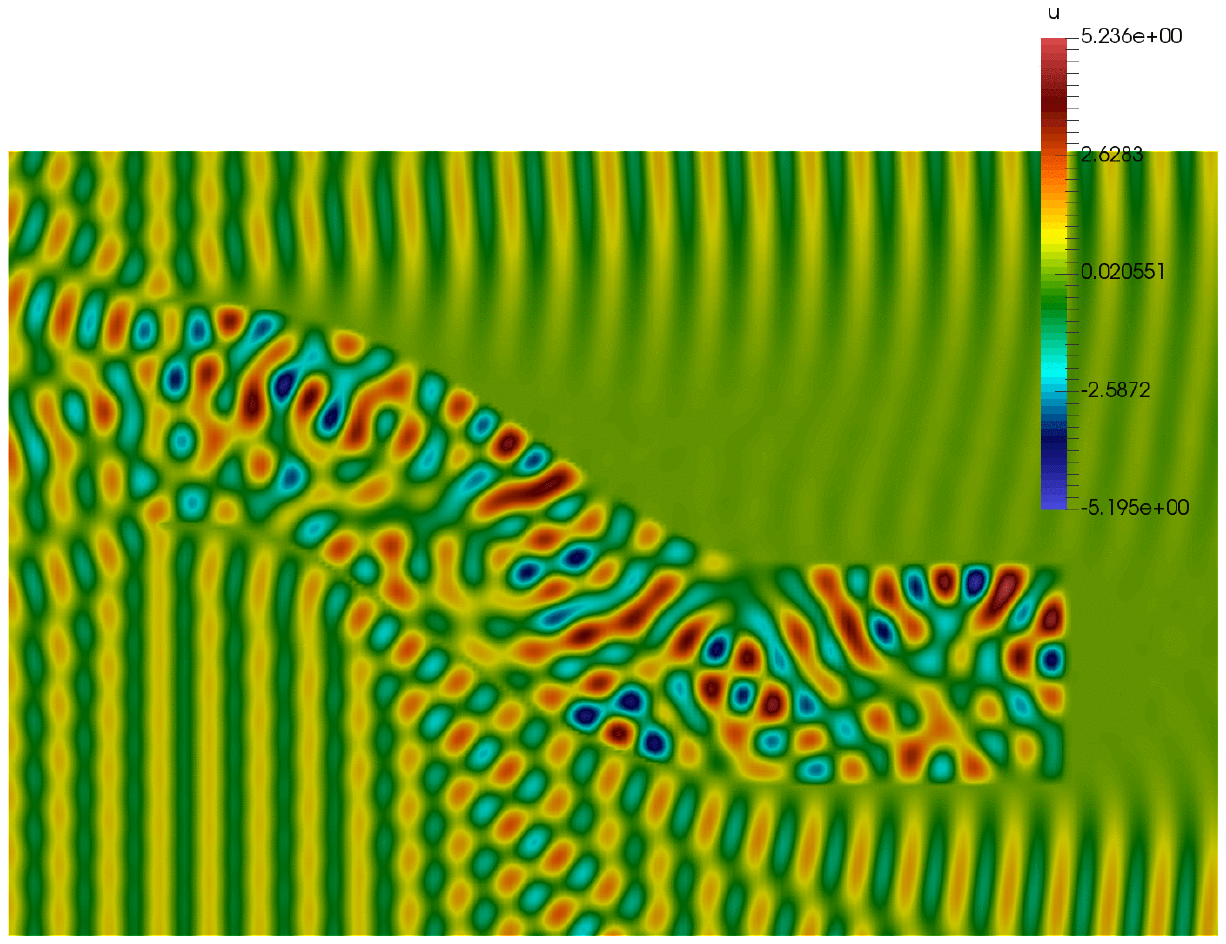}}
	\caption{The COBRA cavity problem.}
	\label{fig:cobra}
\end{figure}

The COBRA cavity problem consists of a plane wave incident upon a curving cavity aperture and scattering in 3D. The geometry is shown within an example solution in \Cref{fig:cobra}. The problem was devised by EADS Aerospatiale Matra Missiles for Workshop EM-JINA 98 and is described in \cite{liu2003scattering,jin2015finite}. The COBRA cavity problem is a benchmark problem in electromagnetism (when we use the time-harmonic Maxwell's equations as the underlying model) but here we present the Helmholtz version. The main difficulty comes from the presence of metallic curved waveguide which can cause trapping in addition to the inherent difficulties present in mid-high frequency regime wave propagation problems. This problem is naturally three-dimensional but a similar two-dimensional test can be designed making use of a cross-section of the cavity.

We first consider discretisation of the problem with 5 points per wavelength using P2 finite elements. In \Cref{Table:Cobra-P2-nppwl5-non-spectral} we display results for the underlying one-level method and the two-level coarse grid approach. The one-level method performs poorly as expected, however, we also see that the coarse grid approach lacks robustness in this setting and performs poorly for larger wave numbers.

\begin{table}[H]
	\centering
	\caption{Results using the one-level and coarse grid methods for the COBRA cavity problem when using 5 points per wavelength, varying the wave number $k$ and the number of subdomains $N$.}
	\label{Table:Cobra-P2-nppwl5-non-spectral}
	\tabulinesep=1.2mm
	\footnotesize
	\begin{tabu}{c|cccc|cccc}
		& \multicolumn{4}{c|}{One-level} & \multicolumn{4}{c}{Coarse grid} \\
		\hline
		$k \ \backslash \ N$ & 20 & 40 & 80 & 160 & 20 & 40 & 80 & 160 \\
		\hline
		100 & 71 & 83 & 104 & $-$ & 21 & 21 & 22 & $-$ \\
		200 & 161 & 203 & 259 & 303 & 87 & 90 & 91 & 91 \\
		300 & 326 & 429 & 574 & 727 & 277 & 309 & 330 & 562 \\
		400 & 326 & 426 & 581 & 705 & 375 & 450 & 463 & 471
	\end{tabu}
\end{table}

Comparing now the spectral coarse space methods, in \Cref{Table:Cobra-P2-nppwl5-spectral-nu} we observe that, in this scenario, the H-GenEO method performs poorly and is often worse than the one-level method, even when taking fairly large numbers of eigenvectors per subdomain. The DtN method can also perform worse than the one-level method in some cases without sufficiently many eigenvectors being taken. However, the DtN approach begins to show relatively low iteration counts once enough eigenvectors are utilised, especially in the many subdomain cases, and can significantly reduce the iteration counts compared to the one-level or coarse grid approaches. Nonetheless, it appears for this 3D problem that the number of eigenvectors required in order to drive down the overall iteration count is rather high.

\begin{table}[H]
	\centering
	\caption{Results using the H-GenEO and DtN spectral coarse space methods for the COBRA cavity problem when using 5 points per wavelength, varying the wave number $k$, the number of subdomains $N$, and the number of eigenvectors used per subdomain $\nu$.}
	\label{Table:Cobra-P2-nppwl5-spectral-nu}
	\tabulinesep=1.2mm
	\footnotesize
	\begin{tabu}{cc|cccc|cccc}
		& & \multicolumn{4}{c|}{H-GenEO} & \multicolumn{4}{c}{DtN} \\
		\hline
		$\nu$ & $k \ \backslash \ N$ & 20 & 40 & 80 & 160 & 20 & 40 & 80 & 160 \\
		\hline
		\multirow{4}{*}{20} & 100 & 66 & 86 & 107 & $-$ & 47 & 58 & 74 & $-$ \\
		& 200 & 165 & 214 & 284 & 366 & 154 & 199 & 250 & 324 \\
		& 300 & 334 & 453 & 611 & 809 & 326 & 429 & 569 & 742 \\
		& 400 & 342 & 440 & 624 & 810 & 337 & 425 & 582 & 727 \\
		\hline
		\multirow{4}{*}{40} & 100 & 65 & 85 & 108 & $-$ & 32 & 34 & 28 & $-$ \\
		& 200 & 164 & 216 & 291 & 377 & 145 & 187 & 233 & 262 \\
		& 300 & 348 & 468 & 638 & 852 & 319 & 424 & 567 & 768 \\
		& 400 & 366 & 500 & 716 & 981 & 337 & 438 & 594 & 778 \\
		\hline
		\multirow{4}{*}{80} & 100 & 64 & 85 & 108 & $-$ & 16 & 12 & 8 & $-$ \\
		& 200 & 155 & 214 & 285 & 382 & 116 & 132 & 99 & 63 \\
		& 300 & 345 & 461 & 636 & 845 & 301 & 399 & 542 & 573 \\
		& 400 & 408 & 558 & 771 & $\times$ & 335 & 453 & 650 & 918 \\
		\hline
		\multirow{4}{*}{160} & 100 & 63 & 85 & 108 & $-$ & 7 & 5 & 5 & $-$ \\
		& 200 & 150 & 210 & 282 & 382 & 46 & 41 & 17 & 11 \\
		& 300 & 324 & 441 & 617 & 830 & 247 & 263 & 249 & 144 \\
		& 400 & 417 & 547 & 749 & $\times$ & 350 & 476 & 576 & 483 \\
		\hline
		\multirow{4}{*}{320} & 100 & $-$ & $-$ & $-$ & $-$ & $-$ & $-$ & $-$ & $-$ \\
		& 200 & 145 & 207 & 281 & 383 & 10 & 10 & 50 & 55 \\
		& 300 & 308 & 422 & 601 & 817 & 100 & 60 & 45 & 29 \\
		& 400 & 401 & 527 & 730 & 994 & 271 & 233 & 147 & 44
	\end{tabu}
\end{table}

We also consider the case of discretising the problem with 10 points per wavelength, again using P2 finite elements. Results in \Cref{Table:Cobra-P2-nppwl10-non-spectral} show that in this case the coarse grid method performs very well and consistently leads to a small number of iterations being required for convergence; this is in stark contrast to the 5 points per wave case. For the spectral coarse spaces, as detailed in \Cref{Table:Cobra-P2-nppwl10-spectral-nu}, the H-GenEO method remains poor, especially for large wave numbers, while the DtN method can perform reasonably well when enough eigenvectors are employed. In particular, the DtN method can require fewer iterations than the coarse grid method, especially for smaller frequencies, but in general struggles to compete in line with the coarse grid approach in terms of robustness.

\begin{table}[H]
\centering
\caption{Results using the one-level and coarse grid methods for the COBRA cavity problem when using 10 points per wavelength, varying the wave number $k$ and the number of subdomains $N$.}
\label{Table:Cobra-P2-nppwl10-non-spectral}
\tabulinesep=1.2mm
\footnotesize
\begin{tabu}{c|cccc|cccc}
	& \multicolumn{4}{c|}{One-level} & \multicolumn{4}{c}{Coarse grid} \\
	\hline
	$k \ \backslash \ N$ & 20 & 40 & 80 & 160 & 20 & 40 & 80 & 160 \\
	\hline
	50 & 27 & 38 & 47 & 52 & 8 & 8 & 8 & 9 \\
	100 & 80 & 103 & 115 & 147 & 11 & 23 & 11 & 11 \\
	150 & 143 & 181 & 235 & 292 & 16 & 16 & 17 & 17 \\
	200 & 192 & 268 & 308 & 427 & 15 & 15 & 16 & 16
\end{tabu}
\end{table}

\begin{table}[H]
\centering
\caption{Results using the H-GenEO and DtN spectral coarse space methods for the COBRA cavity problem when using 10 points per wavelength, varying the wave number $k$, the number of subdomains $N$, and the number of eigenvectors used per subdomain $\nu$.}
\label{Table:Cobra-P2-nppwl10-spectral-nu}
\tabulinesep=1.2mm
\footnotesize
\begin{tabu}{cc|cccc|cccc}
	& & \multicolumn{4}{c|}{H-GenEO} & \multicolumn{4}{c}{DtN} \\
	\hline
	$\nu$ & $k \ \backslash \ N$ & 20 & 40 & 80 & 160 & 20 & 40 & 80 & 160 \\
	\hline
	\multirow{4}{*}{20} & 50 & 26 & 31 & 41 & 46 & 19 & 22 & 20 & 20 \\
	& 100 & 73 & 94 & 111 & 130 & 63 & 77 & 87 & 112 \\
	& 150 & 128 & 171 & 230 & 278 & 120 & 157 & 209 & 244 \\
	& 200 & 186 & 264 & 315 & 427 & 178 & 246 & 283 & 380 \\
	\hline
	\multirow{4}{*}{40} & 50 & 25 & 31 & 40 & 45 & 13 & 12 & 10 & 8 \\
	& 100 & 72 & 89 & 109 & 132 & 51 & 58 & 66 & 57 \\
	& 150 & 129 & 168 & 228 & 279 & 112 & 130 & 181 & 211 \\
	& 200 & 197 & 279 & 328 & 433 & 163 & 225 & 251 & 330 \\
	\hline
	\multirow{4}{*}{80} & 50 & 24 & 31 & 40 & 44 & 6 & 6 & 5 & 4 \\
	& 100 & 68 & 89 & 106 & 132 & 28 & 28 & 19 & 10 \\
	& 150 & 131 & 167 & 224 & 281 & 87 & 101 & 109 & 74 \\
	& 200 & 200 & 270 & 327 & 428 & 142 & 193 & 218 & 255 \\
	\hline
	\multirow{4}{*}{160} & 50 & 24 & 29 & 40 & 44 & 5 & 4 & 4 & 4 \\
	& 100 & 65 & 87 & 109 & 132 & 11 & 9 & 8 & 5 \\
	& 150 & 126 & 164 & 220 & 279 & 50 & 31 & 21 & 14 \\
	& 200 & 192 & 267 & 324 & 430 & 109 & 134 & 83 & 41 \\
	\hline
	\multirow{4}{*}{320} & 50 & 23 & 30 & 40 & 44 & 6 & 6 & 5 & 5 \\
	& 100 & 63 & 85 & 107 & 133 & 6 & 6 & 6 & 5 \\
	& 150 & 119 & 158 & 219 & 279 & 14 & 9 & 8 & 7 \\
	& 200 & 182 & 246 & 311 & 430 & 44 & 28 & 15 & 11
\end{tabu}
\end{table}

\subsubsection{The COBRA cavity in 2D}
\label{sec:CobraCavity2D}

While the COBRA cavity is a 3D problem, we can additionally consider a 2D cross section, as illustrated in \Cref{fig:cobra-2d-cross-section}, and formulate the problem on such a slice. This allows us to consider a cavity problem in 2D and further consider higher wave numbers.

Results for this 2D COBRA cavity problem at 5 points per wavelength are given in \Cref{Table:Cobra2D-P2-nppwl5-non-spectral,Table:Cobra2D-P2-nppwl5-spectral-nu}. Similarly to the full 3D problem, we see that the coarse grid method again lacks robustness in this scenario as the wave number $k$ increases and the H-GenEO is unable to provide any benefit over the one-level method. The DtN approach is able to provide the lowest iteration counts, even with relatively few eigenvectors per subdomain, and provides a very effective method with more eigenvectors.

\begin{table}[H]
	\centering
	\caption{Results using the one-level and coarse grid methods for the 2D COBRA cavity problem when using 5 points per wavelength, varying the wave number $k$ and the number of subdomains $N$.}
	\label{Table:Cobra2D-P2-nppwl5-non-spectral}
	\tabulinesep=1.2mm
	\footnotesize
	\begin{tabu}{c|cccc|cccc}
		& \multicolumn{4}{c|}{One-level} & \multicolumn{4}{c}{Coarse grid} \\
		\hline
		$k \ \backslash \ N$ & 20 & 40 & 80 & 160 & 20 & 40 & 80 & 160 \\
		\hline
		200 & 74 & 99 & 122 & $-$ & 30 & 30 & 30 & $-$ \\
		400 & 103 & 155 & 211 & 279 & 67 & 69 & 70 & 70 \\
		600 & 204 & 303 & 429 & 668 & 145 & 155 & 158 & 160 \\
		800 & 184 & 287 & 392 & 559 & 197 & 214 & 218 & 219 \\
		1000 & 188 & 278 & 400 & 551 & 242 & 265 & 269 & 270
	\end{tabu}
\end{table}

\begin{table}[H]
	\centering
	\caption{Results using the H-GenEO and DtN spectral coarse space methods for the 2D COBRA cavity problem when using 5 points per wavelength, varying the wave number $k$, the number of subdomains $N$, and the number of eigenvectors used per subdomain $\nu$.}
	\label{Table:Cobra2D-P2-nppwl5-spectral-nu}
	\tabulinesep=1.2mm
	\footnotesize
	\begin{tabu}{cc|cccc|cccc}
		& & \multicolumn{4}{c|}{H-GenEO} & \multicolumn{4}{c}{DtN} \\
		\hline
		$\nu$ & $k \ \backslash \ N$ & 20 & 40 & 80 & 160 & 20 & 40 & 80 & 160 \\
		\hline
		\multirow{4}{*}{20} & 200 & 59 & 85 & 122 & $-$ & 24 & 29 & 17 & $-$ \\
		& 400 & 101 & 135 & 202 & 277 & 78 & 87 & 77 & 40 \\
		& 600 & 177 & 271 & 371 & 531 & 161 & 228 & 254 & 221 \\
		& 800 & 192 & 293 & 390 & 517 & 179 & 274 & 341 & 371 \\
		& 1000 & 213 & 302 & 428 & 597 & 203 & 288 & 415 & 521 \\
		\hline
		\multirow{4}{*}{40} & 200 & 54 & 83 & 107 & $-$ & 6 & 7 & 6 & $-$ \\
		& 400 & 94 & 126 & 192 & 271 & 30 & 16 & 9 & 6 \\
		& 600 & 164 & 240 & 348 & 505 & 84 & 82 & 38 & 16 \\
		& 800 & 178 & 259 & 356 & 482 & 129 & 152 & 103 & 34 \\
		& 1000 & 203 & 275 & 394 & 555 & 166 & 198 & 195 & 102 \\
		\hline
		\multirow{4}{*}{80} & 200 & 53 & 76 & 34 & $-$ & 3 & 5 & 3 & $-$ \\
		& 400 & 86 & 120 & 187 & 261 & 7 & 5 & 4 & 3 \\
		& 600 & 154 & 226 & 331 & 490 & 15 & 8 & 5 & 8 \\
		& 800 & 164 & 244 & 338 & 466 & 37 & 22 & 8 & 5 \\
		& 1000 & 184 & 256 & 370 & 529 & 74 & 36 & 14 & 7 \\
		\hline
		\multirow{4}{*}{160} & 200 & 53 & 28 & 1 & $-$ & 2 & 4 & 1 & $-$ \\
		& 400 & 83 & 119 & 183 & 43 & 3 & 3 & 3 & 3 \\
		& 600 & 146 & 218 & 322 & 480 & 4 & 3 & 3 & 6 \\
		& 800 & 154 & 228 & 333 & 458 & 6 & 5 & 3 & 4 \\
		& 1000 & 177 & 243 & 358 & 520 & 9 & 5 & 4 & 3 \\
		\hline
		\multirow{4}{*}{320} & 200 & 38 & 1 & 1 & $-$ & 2 & 1 & 1 & $-$ \\
		& 400 & 81 & 117 & 43 & 1 & 3 & 2 & 3 & 1 \\
		& 600 & 138 & 211 & 315 & 278 & 3 & 3 & 3 & 5 \\
		& 800 & 147 & 221 & 326 & 448 & 3 & 3 & 3 & 3 \\
		& 1000 & 169 & 233 & 346 & 512 & 4 & 3 & 3 & 3
	\end{tabu}
\end{table}

The case of 10 points per wavelength is covered in \Cref{Table:Cobra2D-P2-nppwl10-non-spectral,Table:Cobra2D-P2-nppwl10-spectral-nu} and, again, we broadly see a picture akin to the full 3D problem. In this setting the coarse grid approach provides good performance but the DtN method yields lower iteration counts when sufficiently many eigenvectors are taken, again giving a very effective approach for this problem.

\begin{table}[H]
	\centering
	\caption{Results using the one-level and coarse grid methods for the 2D COBRA cavity problem when using 10 points per wavelength, varying the wave number $k$ and the number of subdomains $N$.}
	\label{Table:Cobra2D-P2-nppwl10-non-spectral}
	\tabulinesep=1.2mm
	\footnotesize
	\begin{tabu}{c|cccc|cccc}
		& \multicolumn{4}{c|}{One-level} & \multicolumn{4}{c}{Coarse grid} \\
		\hline
		$k \ \backslash \ N$ & 20 & 40 & 80 & 160 & 20 & 40 & 80 & 160 \\
		\hline
		200 & 85 & 124 & 171 & 239 & 10 & 11 & 11 & 11 \\
		400 & 137 & 184 & 270 & 356 & 11 & 11 & 12 & 12 \\
		600 & 253 & 329 & 478 & 705 & 20 & 20 & 20 & 21 \\
		800 & 203 & 323 & 448 & 611 & 18 & 18 & 18 & 18 \\
		1000 & 217 & 316 & 470 & 656 & 19 & 19 & 19 & 19
	\end{tabu}
\end{table}

\begin{table}[H]
	\centering
	\caption{Results using the H-GenEO and DtN spectral coarse space methods for the 2D COBRA cavity problem when using 10 points per wavelength, varying the wave number $k$, the number of subdomains $N$, and the number of eigenvectors used per subdomain $\nu$.}
	\label{Table:Cobra2D-P2-nppwl10-spectral-nu}
	\tabulinesep=1.2mm
	\footnotesize
	\begin{tabu}{cc|cccc|cccc}
		& & \multicolumn{4}{c|}{H-GenEO} & \multicolumn{4}{c}{DtN} \\
		\hline
		$\nu$ & $k \ \backslash \ N$ & 20 & 40 & 80 & 160 & 20 & 40 & 80 & 160 \\
		\hline
		\multirow{4}{*}{20} & 200 & 77 & 114 & 159 & 221 & 36 & 42 & 20 & 16 \\
		& 400 & 136 & 195 & 273 & 367 & 111 & 147 & 172 & 123 \\
		& 600 & 241 & 332 & 454 & 658 & 208 & 265 & 346 & 434 \\
		& 800 & 220 & 338 & 485 & 692 & 186 & 282 & 396 & 533 \\
		& 1000 & 242 & 363 & 561 & 787 & 202 & 293 & 419 & 608 \\
		\hline
		\multirow{4}{*}{40} & 200 & 75 & 108 & 153 & 220 & 10 & 7 & 7 & 5 \\
		& 400 & 137 & 187 & 260 & 354 & 66 & 46 & 21 & 13 \\
		& 600 & 236 & 320 & 439 & 632 & 152 & 167 & 128 & 51 \\
		& 800 & 221 & 349 & 470 & 673 & 159 & 236 & 270 & 182 \\
		& 1000 & 254 & 378 & 545 & 785 & 186 & 269 & 365 & 394 \\
		\hline
		\multirow{4}{*}{80} & 200 & 73 & 104 & 152 & 219 & 4 & 4 & 4 & 4 \\
		& 400 & 130 & 174 & 256 & 347 & 11 & 8 & 6 & 5 \\
		& 600 & 220 & 295 & 414 & 590 & 54 & 22 & 11 & 8 \\
		& 800 & 213 & 329 & 454 & 642 & 88 & 74 & 25 & 13 \\
		& 1000 & 256 & 365 & 528 & 727 & 141 & 147 & 80 & 22 \\
		\hline
		\multirow{4}{*}{160} & 200 & 71 & 104 & 148 & 180 & 2 & 3 & 4 & 3 \\
		& 400 & 121 & 165 & 246 & 342 & 5 & 4 & 4 & 3 \\
		& 600 & 210 & 282 & 396 & 577 & 9 & 5 & 5 & 3 \\
		& 800 & 203 & 308 & 435 & 601 & 13 & 8 & 6 & 4 \\
		& 1000 & 245 & 344 & 501 & 688 & 27 & 14 & 8 & 6 \\
		\hline
		\multirow{4}{*}{320} & 200 & 69 & 103 & 137 & 1 & 2 & 3 & 4 & 1 \\
		& 400 & 116 & 158 & 239 & 338 & 3 & 3 & 4 & 3 \\
		& 600 & 197 & 266 & 375 & 568 & 4 & 3 & 3 & 3 \\
		& 800 & 192 & 295 & 412 & 579 & 5 & 4 & 3 & 3 \\
		& 1000 & 225 & 324 & 464 & 653 & 6 & 5 & 4 & 4
	\end{tabu}
\end{table}

\subsection{The Overthrust problem}
\label{sec:OverthrustProblem}

We will now assess our methods on the 3D $20\times 20 \times 4.65\,$km SEG/EAGE
Overthrust model, as introduced and detailed in \cite{aminzadeh1997overthrust}. We use a homogeneous Dirichlet boundary condition at the surface and first-order absorbing boundary conditions along the other five faces of the model. The source is located at $(2.5,2.5,0.58)\,$km and we perform a simulation with P2 finite elements on tetrahedral meshes for $0.5\,$Hz, $1\,$Hz and $2\,$Hz frequencies. For higher frequencies and large-scale computations on the Overthrust problem utilising the grid coarse space see \cite{Dolean:2020:IFD} and \cite{Dolean:2020:LSF}. The Overthrust model can be seen as the three-dimensional counterpart of the Marmousi model and presents the difficulty of the latter (heterogeneities and high frequencies) at a larger scale since we have an extra space dimension.

We will start by testing the different approaches when using 5 points per wavelength, varying the frequency $f$ and the number of subdomains $N$. Results with the one-level and coarse grid methods are shown in \Cref{Table:Overthrust-P2-nppwl5-non-spectral} while results for the H-GenEO and DtN spectral coarse space methods, varying the number of eigenvectors used per subdomain $\nu$, are provided in \Cref{Table:Overthrust-P2-nppwl5-spectral-nu}. We see that in the case of low resolution the grid coarse space proves to be robust. In this scenario the H-GenEO method fails to bring improvement over the one-level method. On the other hand, the DtN approach performs well provided sufficiently many eigenvectors are included and can generally give an improvement over the coarse grid method in terms of reducing the iteration counts.

\begin{figure}[H]
	\centering
	\subfloat[Velocity data\label{fig:overthrust-velocity}]{\includegraphics[width=0.48\textwidth,trim=15cm 0cm 3cm 48cm ,clip]{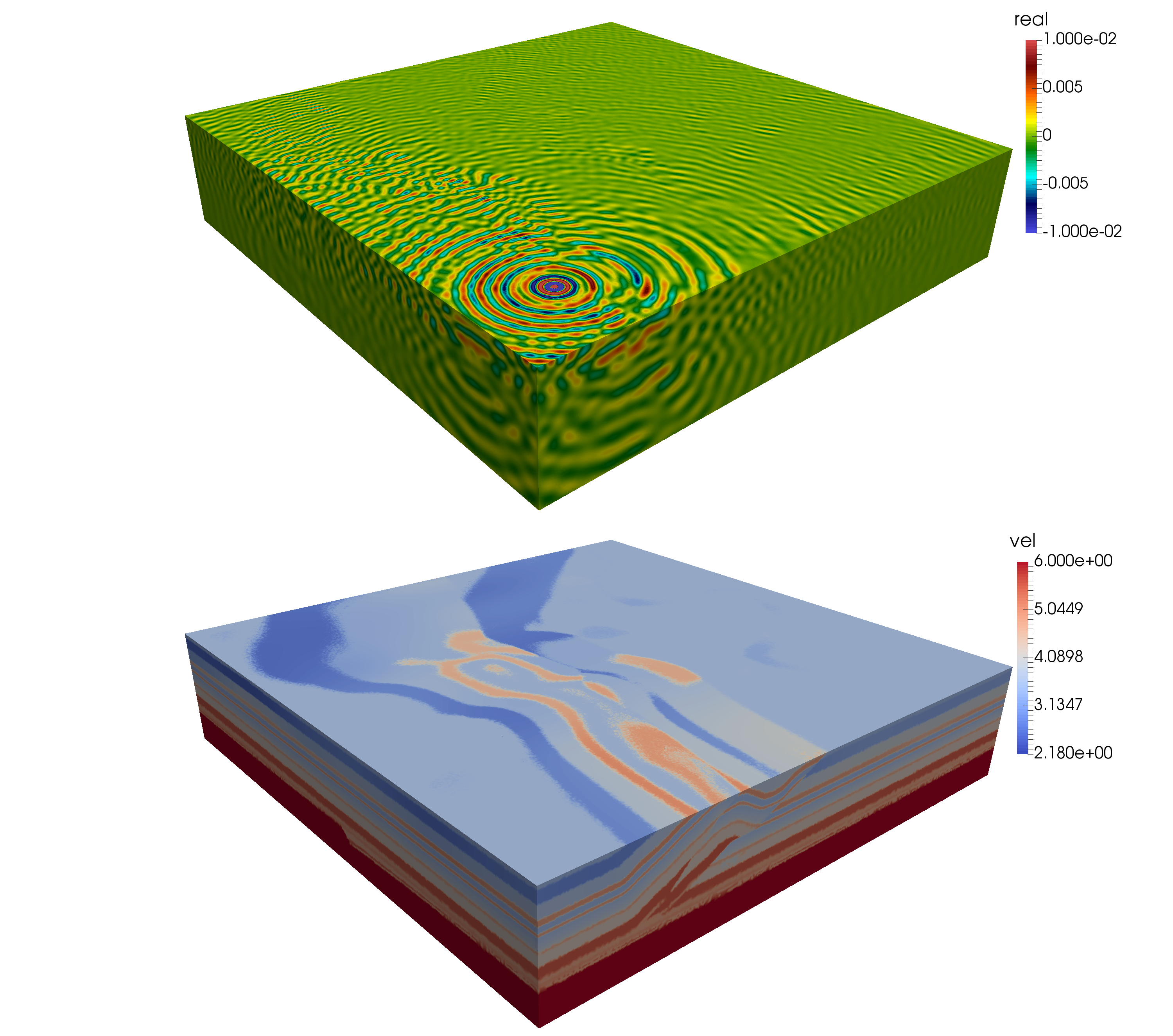}}
	\hfill
	\subfloat[Solution\label{fig:overthrust-solution}]{\includegraphics[width=0.48\textwidth,trim=15cm 46cm 3cm 1cm ,clip]{ovt10adapt004}}
	\caption{The Overthrust problem.}
	\label{fig:overthrust}
\end{figure}

\begin{table}[H]
	\centering
	\caption{Results using the one-level and coarse grid methods for the Overthrust problem when using 5 points per wavelength, varying the frequency $f$ and the number of subdomains $N$.}
	\label{Table:Overthrust-P2-nppwl5-non-spectral}
	\tabulinesep=1.2mm
	\footnotesize
	\begin{tabu}{c|ccccc|ccccc|ccccc}
		& \multicolumn{5}{c|}{One-level (min.\ overlap)} & \multicolumn{5}{c|}{One-level (coarse overlap)} & \multicolumn{5}{c}{Coarse grid} \\
		\hline
		$f \ \backslash \ N$ & 10 & 20 & 40 & 80 & 160 & 10 & 20 & 40 & 80 & 160 & 10 & 20 & 40 & 80 & 160 \\
		\hline
		0.5 & 16 & 21 & 28 & 35 & $-$ & 14 & 17 & 20 & 19 & $-$ & 9 & 10 & 10 & 10 & $-$ \\
		1 & 27 & 36 & 55 & 62 & 78 & 24 & 30 & 48 & 50 & 62 & 17 & 19 & 23 & 23 & 25 \\
		2 & 32 & 47 & 64 & 81 & 104 & 27 & 40 & 53 & 68 & 86 & 18 & 22 & 25 & 27 & 28
	\end{tabu}
\end{table}

\begin{table}[H]
	\centering
	\caption{Results using the H-GenEO and DtN spectral coarse space methods for the Overthrust problem when using 5 points per wavelength, varying the frequency $f$, the number of subdomains $N$, and the number of eigenvectors used per subdomain $\nu$.}
	\label{Table:Overthrust-P2-nppwl5-spectral-nu}
	\tabulinesep=1.2mm
	\footnotesize
	\begin{tabu}{cc|ccccc|ccccc|ccccc|ccccc}
		& & \multicolumn{5}{c|}{H-GenEO (min.\ overlap)} & \multicolumn{5}{c|}{DtN (min.\ overlap)} & \multicolumn{5}{c|}{H-GenEO (coarse overlap)} & \multicolumn{5}{c}{DtN (coarse overlap)} \\
		\hline
		$\nu$ & $f \ \backslash \ N$ & 10 & 20 & 40 & 80 & 160 & 10 & 20 & 40 & 80 & 160 & 10 & 20 & 40 & 80 & 160 & 10 & 20 & 40 & 80 & 160 \\
		\hline
		\multirow{3}{*}{20} & 0.5 & 11 & 14 & 17 & 30 & $-$ & 6 & 6 & 7 & 8 & $-$ & 12 & 15 & 19 & 19 & $-$ & 5 & 5 & 4 & 5 & $-$ \\
		& 1 & 25 & 33 & 47 & 52 & 60 & 21 & 28 & 43 & 45 & 49 & 25 & 34 & 44 & 49 & 56 & 19 & 28 & 42 & 48 & 42 \\
		& 2 & 33 & 47 & 62 & 78 & 99 & 31 & 45 & 59 & 72 & 89 & 32 & 48 & 74 & 101 & 134 & 28 & 41 & 56 & 70 & 97 \\
		\hline
		\multirow{3}{*}{40} & 0.5 & 10 & 12 & 15 & 25 & $-$ & 5 & 5 & 6 & 7 & $-$ & 12 & 15 & 18 & 19 & $-$ & 3 & 4 & 4 & 4 & $-$ \\
		& 1 & 24 & 32 & 44 & 52 & 59 & 17 & 25 & 36 & 38 & 32 & 24 & 34 & 45 & 50 & 56 & 15 & 20 & 28 & 14 & 6 \\
		& 2 & 34 & 50 & 70 & 104 & 123 & 30 & 42 & 55 & 67 & 87 & 35 & 61 & 87 & 112 & 146 & 28 & 42 & 62 & 85 & 138 \\
		\hline
		\multirow{3}{*}{80} & 0.5 & 9 & 11 & 14 & 26 & $-$ & 4 & 5 & 5 & 6 & $-$ & 12 & 14 & 18 & 18 & $-$ & 3 & 3 & 4 & 5 & $-$ \\
		& 1 & 22 & 31 & 44 & 52 & 60 & 13 & 19 & 20 & 14 & 8 & 23 & 33 & 45 & 50 & 56 & 9 & 9 & 10 & 5 & 4 \\
		& 2 & 38 & 61 & 84 & 111 & 139 & 27 & 38 & 53 & 70 & 96 & 41 & 63 & 91 & 116 & 158 & 27 & 44 & 71 & 85 & 65 \\
		\hline
		\multirow{3}{*}{160} & 0.5 & 9 & 11 & 11 & 28 & $-$ & 4 & 5 & 5 & $-$ & $-$ & 12 & 14 & 17 & 19 & $-$ & 4 & 3 & 4 & $-$ & $-$ \\
		& 1 & 19 & 29 & 44 & 52 & 60 & 10 & 10 & 10 & 7 & 7 & 23 & 33 & 46 & 50 & 56 & 7 & 6 & 7 & 5 & 5 \\
		& 2 & 36 & 54 & 77 & 107 & 136 & 23 & 34 & 50 & 56 & 56 & 41 & 62 & 90 & 117 & 159 & 23 & 34 & 42 & 20 & 12
	\end{tabu}
\end{table}

We consider now the case of higher resolution, using 10 points per wavelength with results for smaller frequencies being given in \Cref{Table:Overthrust-P2-nppwl10-non-spectral,Table:Overthrust-P2-nppwl10-spectral-nu}. Here we see that all two-level methods can perform relatively well and both H-GenEO and DtN can yield lower iteration counts than the coarse grid method without needing too many eigenvectors. If we allow for many eigenvectors then the DtN method typically requires the fewest iteration, however, we also note that for fewer eigenvectors the H-GenEO provides a greater reduction in iteration counts than DtN in the $f = 1\,\text{Hz}$ case and, furthermore, relatively little is gained when adding additional eigenvectors.

\begin{table}[H]
	\centering
	\caption{Results using the one-level and coarse grid methods for the Overthrust problem when using 10 points per wavelength, varying the frequency $f$ and the number of subdomains $N$.}
	\label{Table:Overthrust-P2-nppwl10-non-spectral}
	\tabulinesep=1.2mm
	\footnotesize
	\begin{tabu}{c|ccccc|ccccc|ccccc}
		& \multicolumn{5}{c|}{One-level (min.\ overlap)} & \multicolumn{5}{c|}{One-level (coarse overlap)} & \multicolumn{5}{c}{Coarse grid} \\
		\hline
		$f \ \backslash \ N$ & 10 & 20 & 40 & 80 & 160 & 10 & 20 & 40 & 80 & 160 & 10 & 20 & 40 & 80 & 160 \\
		\hline
		0.5 & 19 & 25 & 33 & 39 & 47 & 16 & 21 & 27 & 30 & 36 & 12 & 13 & 14 & 15 & 16 \\
		1 & 31 & 45 & 61 & 75 & 95 & 26 & 39 & 51 & 62 & 78 & 18 & 24 & 25 & 27 & 27
	\end{tabu}
\end{table}

\begin{table}[H]
	\centering
	\caption{Results using the H-GenEO and DtN spectral coarse space methods for the Overthrust problem when using 10 points per wavelength, varying the frequency $f$, the number of subdomains $N$, and the number of eigenvectors used per subdomain $\nu$.}
	\label{Table:Overthrust-P2-nppwl10-spectral-nu}
	\tabulinesep=1.2mm
	\footnotesize
	\hspace*{-0.5cm}
	\begin{tabu}{cc|ccccc|ccccc|ccccc|ccccc}
		& & \multicolumn{5}{c|}{H-GenEO (min.\ overlap)} & \multicolumn{5}{c|}{DtN (min.\ overlap)} & \multicolumn{5}{c|}{H-GenEO (coarse overlap)} & \multicolumn{5}{c}{DtN (coarse overlap)} \\
		\hline
		$\nu$ & $f \ \backslash \ N$ & 10 & 20 & 40 & 80 & 160 & 10 & 20 & 40 & 80 & 160 & 10 & 20 & 40 & 80 & 160 & 10 & 20 & 40 & 80 & 160 \\
		\hline
		\multirow{2}{*}{20} & 0.5 & 8 & 9 & 10 & 11 & 12 & 8 & 7 & 7 & 7 & 7 & 7 & 8 & 10 & 19 & 27 & 7 & 5 & 5 & 5 & 5 \\
		& 1 & 24 & 29 & 34 & 30 & 17 & 23 & 31 & 40 & 44 & 45 & 22 & 30 & 42 & 51 & 64 & 21 & 30 & 38 & 45 & 49 \\
		\hline
		\multirow{2}{*}{40} & 0.5 & 8 & 9 & 9 & 11 & 12 & 6 & 6 & 6 & 6 & 6 & 7 & 8 & 9 & 18 & 26 & 4 & 4 & 4 & 4 & 4 \\
		& 1 & 17 & 18 & 19 & 13 & 14 & 20 & 25 & 29 & 31 & 26 & 19 & 23 & 32 & 44 & 63 & 18 & 25 & 31 & 38 & 43 \\
		\hline
		\multirow{2}{*}{80} & 0.5 & 8 & 8 & 9 & 10 & 12 & 5 & 5 & 5 & 6 & 5 & 7 & 8 & 9 & 18 & 26 & 4 & 3 & 4 & 4 & 4 \\
		& 1 & 12 & 11 & 12 & 11 & 13 & 15 & 18 & 22 & 23 & 15 & 14 & 14 & 21 & 37 & 60 & 15 & 19 & 26 & 24 & 13 \\
		\hline
		\multirow{2}{*}{160} & 0.5 & 7 & 8 & 9 & 10 & 11 & 4 & 5 & 5 & 5 & 7 & 7 & 7 & 8 & 17 & 26 & 3 & 4 & 4 & 4 & 5 \\
		& 1 & 10 & 10 & 10 & 11 & 12 & 11 & 13 & 14 & 17 & 6 & 11 & 12 & 16 & 30 & 57 & 10 & 13 & 12 & 6 & 5
	\end{tabu}
\end{table}

Finally, we examine a higher frequency case of $f = 2\,\text{Hz}$ where the number of degrees of freedom is now $11,\!334,\!869$ and, as such, we consider utilising a larger number of subdomains $N$. Results for the grid coarse space are given in \Cref{Table:Overthrust-P2-nppwl10-non-spectral-2} and show again that the method performs well; for reference, we also report in parentheses the average number of inner iterations required for solving the coarse problem. For the spectral coarse spaces, results are given in \Cref{Table:Overthrust-P2-nppwl10-spectral-nu-2} where we see a similar picture emerge as in the $f = 1\,\text{Hz}$ case. When fewer eigenvectors are used, the H-GenEO approach can converge faster than the DtN method. When more eigenvectors are used, the DtN approach becomes the method that converges fastest and can do so in fewer iterations than the grid coarse space; nonetheless, the H-GenEO method also continues to improve this time and performs similarly to the grid coarse space in terms of iteration counts. Note that the missing entries in \Cref{Table:Overthrust-P2-nppwl10-spectral-nu-2}, for the largest coarse space to be computed, are due to fact that the coarse operator (having dimension $1280\times320$ and dense blocks of size 320) is costly to compute and cannot be handled with current black-box direct solvers.

\begin{table}[H]
	\centering
	\caption{Results using the one-level and coarse grid methods for the Overthrust problem when using 10 points per wavelength for $f = 2\,\text{Hz}$, varying the number of subdomains $N$. The number of degrees of freedom is $11,\!334,\!869$.}
	\label{Table:Overthrust-P2-nppwl10-non-spectral-2}
	\tabulinesep=1.2mm
	\footnotesize
	\begin{tabu}{c|ccc|ccc|ccc}
		& \multicolumn{3}{c|}{One-level (min.\ overlap)} & \multicolumn{3}{c|}{One-level (coarse overlap)} & \multicolumn{3}{c}{Coarse grid} \\
		\hline
		$f \ \backslash \ N$ & 320 & 640 & 1280 & 320 & 640 & 1280  & 320 & 640 & 1280 \\
		\hline
		2 & 149 & 185 & 226 & 125 & 156 & 189 & 24\,(7) & 24\,(9) & 23\,(13)
	\end{tabu}
\end{table}

\begin{table}[H]
	\centering
	\caption{Results using the H-GenEO and DtN spectral coarse space methods for the Overthrust problem when using 10 points per wavelength for $f = 2\,\text{Hz}$, varying the number of subdomains $N$ and the number of eigenvectors used per subdomain $\nu$. The number of degrees of freedom is $11,\!334,\!869$.}
	\label{Table:Overthrust-P2-nppwl10-spectral-nu-2}
	\tabulinesep=1.2mm
	\footnotesize
	\begin{tabu}{c|ccc|ccc|ccc|ccc}
		& \multicolumn{3}{c|}{H-GenEO (min.\ overlap)} & \multicolumn{3}{c|}{DtN (min.\ overlap)} & \multicolumn{3}{c|}{H-GenEO (coarse overlap)} & \multicolumn{3}{c}{DtN (coarse overlap)} \\
		\hline
		$\nu \ \backslash \ N$ & 320 & 640 & 1280 & 320 & 640 & 1280 & 320 & 640 & 1280 & 320 & 640 & 1280 \\
		\hline
		40 & 69 & 73 & 77 & 102 & 97 & 106 & 170 & 209 & 234 & 103 & 146 & 187 \\
		80 & 29 & 39 & 41 & 66 & 70 & 64 & 149 & 206 & 245 & 123 & 124 & 59 \\
        160 & 21 & 28 & 30 & 50 & 61 & 53 & 118 & 189 & 249 & 63 & 17 & 6 \\
        320 & 18 & 26 & $-$ & 42 & 33 & $-$ & 97 & 161 & $-$ & 9 & 7 & $-$ \\
	\end{tabu}
\end{table}

\subsection{Discussion of numerical results}
\label{sec:DiscussionOfNumericalResults}

Numerical assessment of the three coarse spaces considered shows that, depending on the physical parameters (low or high frequency) or on the resolution considered (low resolution at 5 points per wavelength or high resolution at 10 points per wavelength), there is no systematic advantage of one method over another, as shown in the overview in \Cref{Table:Overview}. The underlying problem characteristics (spatial heterogeneities, low or high frequency) can play a strong role in the performance of each of these methods and each method has potential tuning parameters that can affect their efficacy, such as the splitting parameter $s$ for the grid coarse space (which is maintained throughout these tests at its best value, namely $2$, which gives the finest possible coarse space), or the number of local eigenvectors requested for the spectral coarse spaces. The resolution, in turn, can be dictated by the fact that the linear solver is used in the context of an inverse or a direct problem. Further, very often in geophysical computations, especially in full waveform inversion (FWI), one needs to solve repeatedly a linear system with multiple right-hand sides. In these situations, the additional precomputation required by the spectral coarse spaces becomes less of a burden compared to the inner iterations required to solve the coarse problems when using the grid coarse space. Although apparently more expensive as they require solving local eigenvalue problems, spectral methods are potentially very promising in the context of inverse problems.

Nonetheless, we can provide some general observations from our numerical tests. The coarse grid approach is overall fairly robust in terms of the outer iteration counts required and for high wave numbers, so long as enough points per wavelength are used for discretisation. This method tended to be most favourable in the 10 points per wavelength scenarios, which is considered to be too high a level of precision in large-scale geophysical computations as it leads to a substantial, often prohibitive, number of degrees of freedom. The H-GenEO approach, while able to give some positive results for the geophysical test cases at sufficient resolution, failed to be robust for the cavity problems tested, suggesting it may only be a suitable solver method for certain types of problems. Nonetheless, in some scenarios it provides a greater reduction in iteration counts compared to the DtN method when only a relatively small number of eigenvectors are available. The DtN approach performed relatively well in our tests, arguably being the most robust method overall, so long as sufficiently many eigenvectors were utilised and, in particular, tended to be the favourable method in the 5 points per wavelength scenarios. Nonetheless, for the large 3D tests the method started to struggle somewhat unless a large number of eigenvectors were computed, but in this case the approach ultimately becomes rather memory demanding and it is challenging to deal with the coarse operator, which has large dense blocks.

\begin{table}[H]
	\centering
	\caption{An overview of which coarse spaces perform well in the different problem scenarios tests. A \tick~indicates that the method performs well, with \tick\tick~indicating this method was most favourable in a particular instance (note that for the Overthrust problem it is difficult to say that any one method is uniformly the most favourable in each case and so we omit this). A \cross~indicates that a method provided relatively little to no gain over the one-level method.}
	\label{Table:Overview}
	\tabulinesep=1.2mm
	\footnotesize
	\begin{tabu}{c|c|c|c|c|c|c|c|c}
		\multicolumn{3}{c|}{} & \multicolumn{2}{c|}{Coarse grid} & \multicolumn{2}{c|}{H-GenEO} & \multicolumn{2}{c}{DtN} \\
		\hline
		Problem & $d$ & freq & 5 ppwl & 10 ppwl & 5 ppwl & 10 ppwl & 5 ppwl & 10 ppwl \\
		\hline
		\multirow{2}{*}{Marmousi} & \multirow{2}{*}{2D} & low & \tick & \tick & \tick & \tick & \tick\tick & \tick\tick \\
		\cline{3-9}
		& & high & \tick\tick & \tick & \cross & \tick\tick & \tick & \tick \\
		\hline
		\multirow{2}{*}{COBRA Cavity} & \multirow{2}{*}{2D} & low & \tick & \tick & \cross & \cross & \tick\tick & \tick\tick \\
		\cline{3-9}
		& & high & \cross & \tick & \cross & \cross & \tick\tick & \tick\tick \\
		\hline
		\multirow{2}{*}{COBRA Cavity} & \multirow{2}{*}{3D} & low & \tick & \tick\tick & \cross & \cross & \tick\tick & \tick \\
		\cline{3-9}
		& & high & \cross & \tick\tick & \cross & \cross & \tick\tick & \tick \\
		\hline
		\multirow{2}{*}{Overthrust} & \multirow{2}{*}{3D} & low & \tick & \tick & \cross & \tick & \tick & \tick \\
		\cline{3-9}
		& & high & \tick\ & \tick & \cross & \tick & \tick & \tick
	\end{tabu}
\end{table}

In order to give an idea of the advantages and drawbacks of the coarse grid and spectral methods in terms of computational time, we report timings for the setup and solution phases for each method for the Marmousi test case with 10 points per wavelength in~\Cref{Table:timings}. The frequency is $f = 20\,\text{Hz}$ and the number of degrees of freedom is $7,\!813,\!665$. Timings are reported for 40, 80, and 160 subdomains,  with $\nu = 80$ and $\nu = 160$.

\begin{table}[H]
	\centering
	\caption{Reported timings (in seconds) for the Marmousi problem when using 10 points per wavelength at a frequency of $f = 20\,\text{Hz}$ for several choices of the number of subdomains $N$ and number of eigenvectors used per subdomain $\nu$. The number of degrees of freedom is $7,\!813,\!665$. The \emph{setup} column for the one-level method corresponds to the assembly and factorisation of local matrices. The \emph{CS setup} columns correspond to the assembly of the necessary operators for the second level of each two-level method; for the spectral coarse spaces this includes the solution of the local eigenvalue problems as well as the assembly and factorisation of the coarse grid matrix. The \emph{GMRES} columns correspond to the solution phase, with the number of iterations recalled in parenthesis.}
	\label{Table:timings}
	\tabulinesep=1.2mm
	\footnotesize
	\begin{tabu}{c|c|c|cc|c|cc|c|cc|c|cc}
		\multicolumn{2}{c|}{} & \multicolumn{3}{c|}{One-level} & \multicolumn{3}{c|}{Coarse grid} & \multicolumn{3}{c|}{H-GenEO} & \multicolumn{3}{c}{DtN} \\
		\hline
		$N$ & $\nu$ & setup & \multicolumn{2}{c|}{GMRES} & CS setup & \multicolumn{2}{c|}{GMRES} & CS setup & \multicolumn{2}{c|}{GMRES} & CS setup & \multicolumn{2}{c}{GMRES} \\
		\hline
		\multirow{2}{*}{40} & 80 & \multirow{2}{*}{14.9s} & \multirow{2}{*}{113.6s} & \multirow{2}{*}{(198)} & \multirow{2}{*}{1.8s} & \multirow{2}{*}{81.9s} & \multirow{2}{*}{(48)} & 148.8s & 43.0s & (78) & 123.8s & 52.6s & (93) \\
		\cline{2-2} \cline{9-14}
		& 160 & & & & & & & 478.9s & 23.6s & (34) & 377.7s & 31.3s & (47) \\
		\hline
		\multirow{2}{*}{80} & 80 & \multirow{2}{*}{6.9s} & \multirow{2}{*}{86.2s} & \multirow{2}{*}{(272)} & \multirow{2}{*}{0.8s} & \multirow{2}{*}{31.2s} & \multirow{2}{*}{(49)} & 84.2s & 14.4s & (53) & 61.2s & 29.0s & (100) \\
		\cline{2-2} \cline{9-14}
		& 160 & & & & & & & 338.6s & 10.2s & (25) & 189.8s & 21.7s & (56) \\
		\hline
		\multirow{2}{*}{160} & 80 & \multirow{2}{*}{3.2s} & \multirow{2}{*}{72.7s} & \multirow{2}{*}{(389)} & \multirow{2}{*}{0.4s} & \multirow{2}{*}{18.6s} & \multirow{2}{*}{(49)} & 52.5s & 7.0s & (39) & 33.8s & 20.2s & (107) \\
		\cline{2-2} \cline{9-14}
		& 160 & & & & & & & 270.9s & 7.1s & (19) & 124.5s & 20.7s & (59) \\
	\end{tabu}
\end{table}

We see that the setup time of the one-level method, corresponding to the assembly and factorisation of local matrices, decreases accordingly from 14.9s to 3.2s as the size of the subdomains shrinks. The setup time for the coarse space of the coarse grid method follows a similar trend with much smaller timings, from 1.8s to 0.4s, as it corresponds to the assembly and factorisation of the smaller local matrices of the problem defined on the coarse mesh. The GMRES solution time of the one-level method shows poor scaling as the number of subdomains increases, going from 113.6s to 72.7s, because the number of iterations increases from 198 to 389. In contrast, the number of iterations of the coarse grid method is stable as $N$ grows, and the solution time shows good scaling from 81.9s to 18.6s.

For this Marmousi test case, we see that the better spectral coarse space is given by the H-GenEO method: for example, when $\nu = 160$ the number of iterations goes from 34 to 19 as $N$ grows from 40 to 160, whereas for the DtN method it goes from 47 to 59. In terms of GMRES solution time, H-GenEO performs significantly better than the coarse grid method thanks to the significant reduction in the number of iterations and shows decent scaling: from 23.6s to 7.1s when $\nu = 160$. The DtN method performs more poorly and exhibits a worse scaling than the H-GenEO method, primarily because of the increased number of iterations: solution time when $\nu = 160$ goes from 31.3s to 20.7s. However, it performs better than (or, when $N=160$, similar to) the coarse grid method, even though the number of iterations is higher. This can be explained by the relative cost of applying the coarse correction operator $Q$ for the two methods: as mentioned in~\Cref{sec:ImplementationDetails}, the coarse problem in the coarse grid method is solved iteratively with inner GMRES iterations preconditioned by a one-level method, whereas in the spectral methods $E^{-1}$ is applied using its $LU$ factorisation, which is cheaper at these coarse space sizes.

Note also that both spectral methods do not show any gain in solution time at $N = 160$ when $\nu$ goes from 80 to 160. This can be explained by the increased cost of applying $E^{-1}$ as its size grows, relative to the cost of the rest of the iteration (mainly applying the factorisation of the local matrices for the one-level preconditioner, which is cheaper as $N$ grows). That is why there is no decrease in the GMRES solution time even though the number of iterations is nearly halved.

Even though the spectral methods (and especially the H-GenEO method) give good results in terms of number of iterations and solution time, their main drawback is the setup cost of the coarse operator. Indeed, the reported setup costs, covering the solution of the local eigenvalue problems and the assembly and factorisation of the coarse grid matrix, are rather high, especially for large $\nu$: ranging from 478.9s to 270.9s as $N$ grows from 40 to 160 for the H-GenEO method when $\nu = 160$. This can be explained by the significant computational overhead of computing a larger number of eigenvalues and eigenvectors of the local operators, which in practice is far from linear in the number of requested eigenvectors: for example, when $N=160$ we can see the total setup cost going from 52.5s to 270.9s for the H-GenEO method when requesting 160 eigenvectors instead of 80. Interestingly, we can see that the solution of the eigenvalue problems is faster for the DtN method than for the H-GenEO method for the same number of requested eigenvectors, although the resulting coarse space is less efficient in terms of the number of iterations required.

Although the solution of the local eigenvalue problems is completely parallel and shows good scaling as the size of the subdomains shrinks (for example, the setup cost of the H-GenEO method reduces from 148.8s to 52.5s as $N$ grows from 40 to 160 when $\nu = 80$), the high setup cost of the spectral methods can appear prohibitive. However, the setup cost could be leveraged when solving for multiple right-hand sides (for example when computing the radar cross section of an object), or using deflation or recycling techniques when solving successive linear systems that are close to each other.

Deriving full complexity models for our preconditioners is challenging since
many parameters must be taken account, such as the cost of matrix--vector products,
exact subdomain factorisations, eigensolves, and so on. Interested readers are
referred to~\cite{theo} for some bounds on the memory and operation costs of
numerical factorisation and forward plus backward substitutions of MUMPS for
Helmholtz equation. These mostly exhibit that the cost of factorising and
applying the preconditioner will decrease superlinearly for a fixed-size global
problem and an increasing number of subdomains, since the latter will become
smaller and smaller. However, this also makes the coarse operator grow in size.
A multilevel extension thus sounds very appealing in theory. While there are
tools available for symmetric definite problems~\cite{AldGJT19}, they do not
trivially translate to the Helmholtz equation. For eigenanalysis, we rely on
the Krylov--Schur method~\cite{Stewart2002}. There may be more appropriate
mathematical tools for computing a large number of eigenvectors, but this remains
an open problem. In terms of communication between levels for spectral methods,
our implementation uses optimised MPI routines such as MPI\_Gather and
MPI\_Scatter, implementation details may be found
in~\cite{jolivet2013scalable}, for which there are multiple performance
models~\cite{1420226}.

\section{Conclusions}
\label{sec:Conclusions}

We have presented in this work a comparison of several two-level overlapping Schwarz methods on a number of challenging problems of interest from the literature for the Helmholtz problem. Our results illustrate that each of these methods have pros and cons depending on the problem at hand and the particular numerical setting. In particular, no method establishes itself as the superior choice within a wide range of settings, though in particular cases we observe that a certain approach can be more favourable.

Note that our tests on these well-known benchmark problems can be further built upon and the conclusions refined. We have extensively measured timings or considered optimised implementations of the eigenvalue solvers and, in the case of the grid coarse space, the coarse grid problem is the finest possible which further requires another domain decomposition method to solve it as an inner iteration. This study is a first step in the assessment of the state-of-the-art for two-level solvers aimed at time-harmonic wave propagation problems. Additional studies can be directed towards both the larger numerical context (taking into account higher approximation order, for example) and further applicative contexts (the solution of inverse problems in FWI, or the time-harmonic model of elastic waves which is at the same time more complex and more physically accurate). We remark that this is predominantly a numerical study and, because of the mathematical difficulty of the models involved, a theoretical insight appears to remain out of reach, with a particular challenge being to cover the complexity of all practical situations of interest. Nonetheless, a general conclusion to be drawn is that coarse spaces are clearly a key feature in achieving scalability and robustness and, further, that the heuristics required for time-harmonic wave propagation problems are necessarily very different from the established methodology in the case of symmetric and positive definite problems.

\bibliographystyle{abbrv}
\bibliography{preprint}

\end{document}